\documentclass[11pt,a4paper,oneside]{article}

\usepackage{graphicx}
\usepackage[T1]{fontenc}
\usepackage[english]{babel}
\usepackage{url}
\usepackage{textcomp}
\usepackage[noBBpl]{mathpazo}
\usepackage{booktabs}
\usepackage{amsmath}
\usepackage{amsfonts}
\usepackage{amsthm}

\usepackage{natbib}
\usepackage{hypernat}
\usepackage[dvips,ps2pdf,pagebackref=false,colorlinks=true,citecolor=black,linkcolor=black,plainpages=false,pdfpagelabels=true]{hyperref}
\usepackage[all]{hypcap}

\graphicspath{{.}{./figures/}}

\newcommand{\thefont}[2]{\fontsize{#1}{#2}\fontshape{n}\selectfont}
\newcommand{\1}{\rlap{\thefont{10pt}{12pt}1}\kern.16em\rlap{\thefont{11pt}{13.2pt}1}\kern.4em}

\renewcommand{\baselinestretch}{1.5}

\makeatletter
\def\singlespace{\def\baselinestretch{1}\@normalsize}

\newcommand{\argmax}{\displaystyle \mathop{argmax}}

\renewcommand{\epsilon}{\varepsilon} 
\renewcommand{\hat}{\widehat} 

\def\1{\mbox{1\hspace{-.35em}1}}
\def\R{\mathbb{R}}

\def\P{\mathbb{P}}
\def\E{\mathbb{E}}

\def\R{\mathbb{R}}

\def\Z{\mathbb{Z}}

\newcommand{\be}{\mbox{${\mathbf e}$}}
\newcommand{\bz}{\mbox{${\mathbf z}$}}

\newcommand{\bff}{\mbox{${\mathbf F}$}}

\newcommand{\bX}{\mbox{${\mathbf X}$}}

\newcommand{\bet}{\mbox{\boldmath$\eta$}}

\newcommand{\bbeta}{\mbox{\boldmath$\beta$}}
\newcommand{\btheta}{\mbox{\boldmath$\theta$}}
\newcommand{\bxi}{\mbox{\boldmath$\xi$}}

\newcommand{\la}{\mbox{$\langle$}}
\newcommand{\ra}{\mbox{$\rangle$}}

\newtheorem{theorem}{Theorem}
\newtheorem{Proposition}{Proposition}[]
\newtheorem{corollary}[Proposition]{Corollary}

\title{{\sc Wavelet penalized likelihood estimation in generalized functional models}}

\author{Ir\`ene Gannaz}

\date{2012-10}

\newcounter{comptefigure}

\begin{document}

\setcounter{comptefigure}{0}

\setlength{\parindent}{0pt}
\setlength{\parskip}{\baselineskip}

\begin{center}

{\Large \bf Wavelet penalized likelihood estimation in generalized functional models}

{\bf Ir\`ene Gannaz}

Universit\'e de Lyon\\ \vspace{-0.2\baselineskip}
CNRS UMR 5208\\ \vspace{-0.2\baselineskip}
INSA de Lyon\\ \vspace{-0.2\baselineskip}
Institut Camille Jordan\\ \vspace{-0.2\baselineskip}
20, avenue Albert Einstein\\ \vspace{-0.2\baselineskip}
69621 Villeurbanne Cedex\\ \vspace{-0.2\baselineskip}
France.\\ \vspace{-0.2\baselineskip}
\texttt{irene.gannaz@insa-lyon.fr}\\

\end{center}

\begin{abstract}
The paper deals with generalized functional regression. The aim is to estimate the influence of covariates on observations, drawn from an exponential distribution. The link considered has a semiparametric expression: if we are interested in a functional influence of some covariates, we authorize others to be modeled linearly. We thus consider a generalized partially linear regression model with unknown regression coefficients and an unknown nonparametric function. We present a maximum penalized likelihood procedure to estimate the components of the model introducing penalty based wavelet estimators. Asymptotic rates of the estimates of both the parametric and the nonparametric part of the model are given and quasi-minimax optimality is obtained under usual conditions in literature. We establish in particular that the $\ell^1$-penalty leads to an adaptive estimation with respect to the regularity of the estimated function. An algorithm based on {\it backfitting} and Fisher-scoring is also proposed for implementation. Simulations are used to illustrate the finite sample behaviour, including a comparison with kernel and splines based methods.

\end{abstract}



\medskip

\section{Introduction} 

  The aim of this paper is to consider a regression model, where the response $Y$ is to be predicted by covariates $\bz$, with $Y$ real-valued and with $\bz$ a real explanatory vector. Relaxing the usual assumption of normality, we consider a generalized framework. The response value $y$ is drawn from a one-parameter exponential family of distributions, with a probabilistic density of the form: \begin{equation}\label{modele}\exp\left(\frac{y\eta(\bz)-b(\eta(\bz))}{\phi}+c(y,\phi)\right).\end{equation} In this expression, $b(\cdot)$ and $c(\cdot)$ are known functions, which determine the specific form of the distribution. The parameter $\phi$ is a dispersion parameter and is also supposed to be known in what follows. The unknown function $\eta(\cdot)$ is the natural parameter of the exponential family, which carries information from the explanatory variables. Given a random sample of size $n$ drawn independently from a generalized regression model, the aim is to predict the function $\eta(\cdot)$. Such a model gives a large scope of applications because observations can result from many distribution families such as Gaussian, Poisson, Binomial, Gamma, {\it etc}. For a more thorough description of generalized regression modeling, we refer to \cite{MacCullaghNelder} or \cite{FahrmeirTutz}.

In the following $\bz=(\bX,t)$, with $\bX$ a $p$-dimensional vector and $t$ a real-valued covariate. The function $\eta(\cdot)$ is given by:
  \begin{equation}\label{GPLM}
   G(\mu(\bX,t))=\eta(\bX,t)=\bX^T\bbeta+f(t),
  \end{equation}
where $\bbeta$ is an unknown $p$-dimensional real parameter vector and $f(\cdot)$ is an unknown real-valued function; such a model is called a generalized partially linear model (GPLM). Taking $\bX^T=0$ or $f=0$ lead respectively to a functional and a linear model. Results presented here can thus easily be deduced for both models. 
Given the observed data $(Y_i,\bX_i,t_i)_{i=1,\ldots,n}$, the aim is to estimate from the data the vector $\bbeta$ and the function $f(\cdot)$.

  In a Gaussian modelization, \cite{Rice} and \cite{Speckman} put in evidence that the rates of the estimates for linear and nonlinear parts could not be both optimized without a control of the correlation between the explanatory variable of the linear part and the functional part of the model. With such a control, \cite{Speckman} proves that it is possible to obtain both optimal linear rate and nonparametric rate for the estimates. To my knowledge, the only paper establishing such a result in GPLM is \cite{MammenVanDerGeer}.

  Many papers focus on the asymptotic behaviour of the estimator of the parametric part $\bbeta$ in generalized partially linear models (see {\it e.g.}~\cite{Chen87} by a penalized {\it quasi-least squares} or \cite{SeveriniStaniswalis} by {\it profile likelihood} methods). A recent article of \cite{BoenteHeZou} establishes a uniformly convergent estimation for $\bbeta$ using a robust {\it profile likelihood} similar to \cite{SeveriniStaniswalis}. But few works consider simultaneously the parametric and the nonparametric part of the model. The paper of \cite{MammenVanDerGeer} shows minimax optimality for the estimations of both $f$ and $\bbeta$ with a penalized {\it quasi-least squares} procedure with a Sobolev type penalty. Their correlation condition given for attaining optimality for both parametric and nonparametric estimators appear to be less restrictive than Rice's (1986) or Speckman's (1986).

  This paper proposes a new estimation procedure based on wavelet decompositions. Wavelet based estimators allow to consider heterogeneous functions and can lead to adaptive procedure with respect to the smoothness of the target function. The use of wavelets for estimating the nonparametric component of a Gaussian partially linear model has been investigated by \cite{Meyer03}, \cite{ChangQu}, \cite{FadiliBullmore} or \cite{gannaz} and \cite{gannaz_these} more recently. But it has not been studied in the context of GPLM. Estimation methods encountered in literature need the choice of a smoothing parameter, whose optimal value depends on the regularity of the functional part $f$. A cross-validation procedure is then necessary to evaluate this parameter. As noted by \cite{Rice} or \cite{Speckman}, cross-validation can present much instability in partially linear models. The use of wavelets here leads to a procedure where no cross-validation is needed. This adaptivity is the main novelty of our estimation scheme in such models. We moreover establish the near-minimax optimality of the estimation for both the linear predictor $\bbeta$ and the nonparametric part $f$, under usual assumptions of correlation between the two parts. 
Finally, we present an algorithm for computing the estimates. 
  
 The paper is organized as follows: Section 2 presents the assumptions and the estimation procedure. It also gives the main properties of our estimators. We distinguish two cases: non adaptive penalties and a $\ell^1$ type penalty, which leads to adaptivity. In Section 3, we propose a computational algorithm of the adaptive estimation procedure. We present a simulation study with a numerical comparison with kernel and splines estimators. Proofs of our results are given in the Appendix.

\section{Assumptions and estimation scheme}

  We consider a generalized regression model, where the response $Y$ depends on covariates $(\bX,t)$, where $Y$ is real-valued, $\bX$ is a $p$-dimensional vector and $t$ is a real-valued covariate. The value $y$ is drawn from an exponential family of distributions, with a probabilistic density of the form given by equation (\ref{modele}). As noted above, we are interested in this paper by a semiparametric expression (\ref{GPLM}) of the function $\eta(\cdot)$.
  
The vector $\bbeta$ and the function $f$ are respectively the parametric and nonparametric components of the generalized partially linear model (GPLM). In the following, $(Y_i,\bX_i,t_i)_{i=1,\ldots,n}$ will denote an independent random sample drawn from the GPLM described here.

Recall the conditional mean and variance of the $\text{i}^{\text{th}}$ response $Y_i$ are given by: 
\begin{eqnarray}
\E[Y_i|\bz_i]&=&\dot b(\eta(\bz_i))\,=\,\mu(\bz_i),\\
Var[Y_i|\bz_i]&=&\phi\,\ddot b(\eta(\bz_i)),
\end{eqnarray} where $\dot b(\cdot)$ and $\ddot b(\cdot)$ denote respectively the first and second derivatives of $b(\cdot)$. The function $G=\dot b^{-1}$ is called link function and one has $G(\E(Y_i|\bz_i))=\eta(\bz_i)$.

  \subsection{Penalized maximum likelihood}

The aim of the paper is to estimate simultaneously the parameter $\bbeta$ and the function $f$, given the observed data $(Y_i,\bX_i,t_i)_{i=1,\ldots,n}$. We propose a penalized maximum loglikelihood estimation. Let $\mathcal L$ denotes the loglikelihood function: $\mathcal L(y,\eta)=y\eta-b(\eta).$ We consider throughout the paper that the estimators $\hat f_n$ and $\hat \bbeta_n$ are solutions of:
\begin{equation}\label{critere}
(\hat f_n,\hat\bbeta_n)=\argmax_{\{f,\,\|f\|_\infty\leq C_\infty\},\bbeta\in\R^p} K_{n,\lambda}(f,\bbeta) \end{equation}
$$\text{ with~~} K_{n,\lambda}(f,\bbeta) = \sum_{i=1}^n \mathcal L\left(y_i, \bX_i^T\bbeta+f(t_i)\right) \, - \, \lambda\,Pen(f).$$
{The presence of the penalty $Pen(\cdot)$ corresponds to a constraint on $f$ being in a given space. The parameter $\lambda$ controls the degree of smoothness given on the estimation $\hat f_n$.}

In what follows, we will assume the penalty $Pen(\cdot)$ is convex and the likelihood function is bounded in order to ensure the existence of maxima (unicity is not garanteed but there is no local maxima). We refer to \cite{AnestisGijbelsNikolova} for the conditions of existence of such a maximization problem with non convex penalties. We did not succeed in getting rid of the constraint $\|f\|_\infty\leq C_\infty$ in the proofs but such a condition does not seem too restrictive in practice.

  The computation of the maximization problem is done in two steps. Some studies, such as \cite{Speckman}, incite to estimate first the functional part. Thus, we will proceed as follows:
 \begin{enumerate}
\item $\widetilde f_{n,\bbeta}=\argmax_{\{f,\,\|f\|_\infty\leq C_\infty\}} K_{n,\lambda}(f,\bbeta)$.\\

\item $\hat\bbeta_n=\argmax_{\bbeta\in\R^p} K_{n,\lambda}(\widetilde f_{n,\bbeta},\bbeta)$.\\ 
\item $\hat f_n=\widetilde f_{n,\hat\bbeta_n}$.
\end{enumerate}

Actually, a classical procedure during the second step is to maximise a modified criterion called {\it profile likelihood} (see among others \cite{SeveriniWong} or \cite{BoenteHeZou}). The criterion maximized is then $\sum_{i=1}^n \mathcal L(y_i,\dot b(\bX_i^T\bbeta+\widetilde f_{n,\bbeta}(t_i)))$. The expression of $\dot b(\bX_i^T\bbeta+\widetilde f_{n,\bbeta}(t_i))$ can indeed be simplified using first order conditions of step 1. Due to the non-linearity of our procedure, we choose here to keep a loglikelihood approach.

  Note also that an usual estimation procedure used in generalized models is quasi-likelihood maximization, which has the advantage of being efficient for a wider class of models than generalized. For details, we refer among others to \cite{SeveriniStaniswalis} or \cite{MacCullaghNelder}. In GPLM, quasi-likelihood estimation was developed for example by \cite{Chen87} and \cite{MammenVanDerGeer}.

\subsection{Discrete wavelet transform}

  The aim of the present work is to introduce a wavelet penalty in estimation, which will lead to a wavelet representation of the functional part. Thanks to their time and frequency localization, wavelets allow to capture local singularities and to consider more spatially inhomogeneous functions than kernel or splines based estimations. Moreover some non linear wavelets procedures are adaptive with respect to the smoothness of the estimated function while splines or kernel estimators need to choose a smoothing parameter using a cross-validation procedure. We recall here briefly some facts on wavelets bases. For more precision on wavelets, the reader is referred to \cite{Daubechies}, \cite{Meyer92} or \cite{Mallat99}.

  Let $\left(L^2[0,1],\la \cdot, \cdot \ra\right)$ be the space of squared-integrable functions on $[0,1]$ endowed with the inner product $\la f,g \ra=\int_{[0,1]}f(t)g(t)\,dt$. Throughout the paper we assume that we are working within an $R$-regular ($R\geq 0$) multiresolution analysis of $\left(L^2[0,1],\la \cdot, \cdot \ra\right)$, associated with an orthonormal basis generated by dilatations and translations of a compactly supported scaling function, $\varphi(t)$, and a compactly supported mother wavelet, $\psi(t)$. For simplicity reasons, we will consider periodic wavelet bases on $[0,1]$.
  
  For any $j\geq 0$ and $k=0,1,\ldots,2^j-1$, let us define $\varphi_{j,k}(t)=2^{j/2}\varphi(2^j t-k)$ and $\psi_{j,k}(t)=2^{j/2}\psi(2^j t-k)$. Then for any given resolution level $j_0\geq 0$ the family
$$ \left\{\varphi_{j_0,k},\,k=0,1,\ldots,2^{j_0}-1;\;\psi_{j,k},\,j\geq j_0;\,k=0,1,\ldots,2^j-1\right\} $$ is an orthonormal basis of $L^2[0,1]$. Let $f$ be a function of $L^2[0,1]$ ; if we denote by $c_{j_0,k}=\la f,\varphi_{j_0,k}\ra$ ($k=0,1,\ldots,2^{j_0}-1$) the scaling coefficients and by $d_{j,k}=\la f,\psi_{j,k}\ra$ ($j\geq j_0;\,k=0,1,\ldots,2^j-1$) the wavelet coefficients of $f$, the function $f$ can then be decomposed as follows:
$$ f(t)=\sum_{k=0}^{2^{j_0}-1}c_{j_0,k}\varphi_{j_0,k}(t)+\sum_{j=j_0}^\infty\sum_{k=0}^{2^j-1}d_{j,k}\psi_{j,k}(t),\quad t\in [0,1].$$
  Yet in practice, we are more concerned with discrete observation samples rather than continuous. Consequently we are more interested by the discrete wavelet transform (DWT). Given a vector of real values $\be=(e_1,\ldots,e_n)^T$, the discrete wavelet transform of $\be$ is given by $\btheta=\Psi_{n\times n}\be$, where $\btheta$ is an $n\times 1$ vector comprising both discrete scaling coefficients, $\theta_{j_0,k}^S$, and discrete wavelet coefficients, $\theta_{j,k}^W$. The matrix $\Psi_{n\times n}$ is an orthogonal $n\times n$ matrix associated with the orthonormal periodic wavelet basis chosen, where one can distinguish the {\it Blocs} spanned respectively by the scaling functions and the wavelets.

  Note that if $\bff$ is a vector of function values $\bff=(f(t_1),\ldots,f(t_n))^T$ at equally spaced points $t_i$, then the corresponding empirical coefficients $\theta_{j_0,k}^S$ and $\theta_{j,k}^W$ are related to their continuous counterparts $c_{j_0,k}$ and $d_{j,k}$ with a factor $n^{-1/2}$. Because of the orthogonality of $\Psi_{n\times n}$, the inverse DWT is simply given by $\bff=\Psi_{n\times n}^T\btheta$. If $n=2^J$ for some positive integer $J$, \cite{Mallat} propose a fast algorithm, that requires only order $n$ operations, to compute the DWT and the inverse DWT.

\subsection{Assumptions and asymptotic minimaxity}

 To ameliorate the comprehension of the results and the assumptions, the subscript $0$ will identify in the following the true values of the model. Let $\|.\|$ denotes the euclidean norm on $\R^p$ and $\|h\|_n^2=\frac 1 n \sum_{i=1}^n h(t_i)^2$ for any function $h$. The relevant inner product associated to the exponential law of the data is the following: $$\text{for all~} x\in\R^n,\, y\in\R^n\, <x,y>_G=\frac 1 n \sum_{i=1}^n \ddot b(\bX_i^T\beta_{0}+f_0(t_i))x_i^T y_i.$$ The norm will be noted $\|.\|_G$.
Recall that the function $\ddot b(\cdot)$ is associated to the variance of the observations and that one has $\ddot b>0$.

  Due to the use of the Discrete Wavelet Transform described above, we will consider in the following that the functional part is observed on an equidistant sample $t_i=\frac i n$, and that the sample size satisfies $n=2^J$ for some positive integer $J$. Even if it is quite restrictive in practice, many applications verify such assumptions. For instance in neurosciences, electrocadiogram studies, functional MRI, {\it etc}, the signals are observed at regular time intervals and $t_i$ denotes the time.

We introduce the notation $\eta_{0,i}=\bX_i^T\bbeta_{0}+f_0(t_i)$ for $i=1,\dots,n$ and $\bet_0=(\eta_{0,i})_{i=1,\dots,n}$. Let $B$ be the diagonal matrix with $i^{th}$ term $\ddot b(\eta_{0,i})$. Define $H=X(X^TBX)^{-1}X^T$ the projection matrix for the inner product $<.,.>_G$ on the space generated by the columns of $X$. The hat matrix $H$ admits a rank and thus a trace equal to $p$. If $h_i=\bX_i(X^TBX)^{-1}\bX_i^T$ denotes the $i^{th}$ diagonal term of $H$, this means that $\sum h_i=p$.

We introduce the following assumptions: 
\begin{description}
\item[(A1)] $\frac{1}{n} \sum_{i=1}^n\ddot b(\eta_{0,i}) \bX_i \bX_i^T$ converges to a strictly positive matrix when $n$ goes to infinity, and $\frac 1 n X^T F_0$ goes to 0 when $n$ goes to infinity, with $F_0=\left(f_0(t_1),\ldots,f_0(t_n)\right)^T$.

\vspace{10pt}

\item[(A2)] $h=\displaystyle{\max_{i=1\dots n}}\, h_i\to 0$.

\vspace{10pt}

\item[(A3)]  $\sup_{\{\bet\in\R^n,\,\|\bet-\bet_0\|_n\leq 2 C_\infty\}}\displaystyle{\sup_{i=1,\dots,n}} \dddot b(\eta_{i})\,\leq\,\dddot b_{\infty}\,<\,\infty$. 

\vspace{10pt}

\item[(A4.1)] There exists a constant $a>0$ such that 
 $\displaystyle{\max_{i=1,\dots,n}}\E\left[\exp(\dot{\mathcal L} (Y_i,\eta_{0,i})^2/a)\right]\leq a,$

\vspace{10pt}
 
\item[(A4.2)] There exist constants $K,\,\sigma_0^2\,>0$ such that 
$\displaystyle{\max_{i=1,\dots,n}} K^2\left(\E[\exp\left(\left|\dot{\mathcal L}(Y_i,\eta_{0,i})\right|/K^2\right)]-1\right)\leq \sigma_0^2.$
\end{description}

Assumption (A1) ensures the identifiability of the model. 
Assumption (A3) is not restrictive and holds for classical distributions like Gaussian, Binomial, Poisson, {\it etc}. Some more restrictive assumptions are made on the form of the distribution with assumptions (A4.1) and (A4.2). Assumption (A4.1) corresponds to exponential tails and is weaker than assumption (A4.2), which corresponds to sub-Gaussian tails. When assumptions (A4.1) or (A4.2) hold, they imply that there exists a positive constant $\ddot b_\infty$ such that $\max_i \ddot b(\bX_i^T\beta_0+f_0(t_i))\,\leq\,\ddot b_\infty\,<\,\infty$. As a consequence, for every vector $v\in\R^n$, $\|v\|_G\leq \ddot b_\infty \|v\|_n$.

In \cite{MammenVanDerGeer}, the authors establish their asymptotic results for distributions with a bounded likelihood and for the Binomial distribution. Comparing to their conditions, assumptions on distributions are here less restrictive. For example, they hold for the Poisson model, which does not verify \cite{MammenVanDerGeer}'s assumptions.

  We then aim to control the correlation between the linear and the nonparametric parts of the model. Following \cite{Rice} or \cite{Speckman} we decompose the components of the design matrix $X$ into a sum of a deterministic function of $L^2[0,1]$ and a noise term. More precisely, the $(i,j)$-component of $X$, say $x_{i,j}$, is supposed to take the form $x_{i,j}=g_j(t_i)+\xi_{i,j}$ with functions $(g_j)_{j=1,\ldots,p}$ forming an orthogonal family on $\left(L^2[0,1],\la\cdot,\cdot\ra\right)$ and with $\xi_{i,j}$ denoting a realization of a random variable $\bxi_j$. 

  We make an assumption on the distribution of the random variables $\bxi_j$, and of course we suppose we control the regularity of the functions $g_j$. 
\begin{description}
\item[($A_{corr}$)] $\forall j=1,\dots,p,\,i=1,\dots,n,\,\,X_{i,j}=g_j(t_i)+\xi_{i,j}$, with polynomial functions $g_j$ of degree less or equal than the number of vanishing moments of the wavelet considered. For all  $j=1,\dots,p$, $(\xi_{i,j})_{i=1,\dots,n}$ is a $n$-sample such that $\max_{i=1,\dots,n}\E\left[\exp(\bxi_{i,j}^2/a_j)\right]\leq a_j$, 
for given constants $a_j>0$.
\end{description}

This condition corresponds to covariates $x_{i,j}$ and $t_i$ which are drawn from correlated random variables $X_j$ and $T$ with the relation $\E(X_j|T)=g_j(T)+\xi_j$. If $t_i$ denotes the time, this condition allows the covariate $(\bX_{i,j})_{i=1,\dots,n}$ to have a polynomial trend with respect to the time.

The form is very similar to the one given by \cite{Rice}. In the gaussian partially linear regression the author proved that estimators of both the linear part and the functional part could attain optimal rates, respectively for each part, when $f$ belongs to a Sobolev space $\mathcal W^m$, if functions $g_j$ were polynomial with a degree strictly smaller than $m$. The condition $(A_{corr})$ on design covariates appears to be more flexible than Rice's (1986) in the sense that the maximal degree of the polynomial functions intervening in the covariates depends on the number of vanishing  moments of the wavelet, instead of depending on the regularity of the function $f$.

Yet, in the Gaussian framework, the condition of \cite{Rice} has been relaxed afterwards, first by \cite{Speckman} who only suppose that the functions $g_j$ measuring the correlation are $m$ times continuously differentiable. Recently, \cite{diff_finies} also proposed a difference based approach leading to minimax optimality as soon as either $f$ or functions $g_j$ are smooth enough ; more precisely if $g_j$ are $a_g$-Lipschitz and $f$ is $a_f$-Lipschitz than it is sufficient that $a_g+a_f>1/2$. Condition $(A_{corr})$ is more retrictive than those of \cite{Speckman} or \cite{diff_finies} because it imposes the polynomial form on the correlation functions. In the generalized framework, $(A_{corr})$ is also more restrictive than the correlation assumption of \cite{MammenVanDerGeer} where, for example, minimax asymptotics are available if the function $f$ and the functions $g_j$ all belong to a Sobolev space $\mathcal W^m$. Yet, when considering non polynomial correlation functions $g_j$ in simulations, we still obtain a similar behaviour of our estimators. Thus assumption $(A_{corr})$ could probably bee relaxed.

\subsubsection{Nonadaptative case} \label{sec:nonadapt}

We suppose that the function $f_0$ in the nonparametric part belongs to a function class $\mathcal A$ which can be described by $\mathcal A=\left\{g,\,J(g)\leq C\right\}$, with $J(\cdot)$ a given criterion on the functions from $[0,1]$ to the positive real line.  
We suppose that the $\delta$-entropies of the subspace $\mathcal A$ for the distance associated to the norm $\|\cdot\|_n$ behave like $$\lim\sup_{n\to\infty}\,\sup_{\delta>0}\;\delta^\nu\,\mathcal H(\delta,\mathcal A,\|\cdot\|_n)\;<\;\infty,$$ for a given $0<\nu<2$.

The penalty in equation~(\ref{critere}) is chosen according to the following assumptions:
\begin{description}
\item[(A5)] For any function $g$, $\frac{\lambda v_n^2}{n}Pen(g)\geq J(g)$ with $v_n=n^{1/(2+\nu)}$.
\vspace{10pt}
\item[(A6)] $f\mapsto K_{n,\lambda}(f,\bbeta)= \sum_{i=1}^n \mathcal L(y_i,\bX_i^T\bbeta+f(t_i))-\lambda Pen(f)$ is concave. 
\vspace{10pt}
\item[(A7)] The penalty $Pen(h)$ applies only to the wavelet decomposition coefficients $(\theta^W_{j,k})_{j\geq j_S, \,k\in\Z}$ of the function~$h$ for some $0\leq j_S<J-1$.
\end{description}
Note assumption (A7) is introduced to exploit the special structure given in $(A_{corr})$ through a penalty on wavelet coefficients.

We are now in position to give our first asymptotic result.

\begin{theorem}\label{th1}
Suppose assumptions (A1) to (A3), (A4.1), (A5) and (A6) hold and $\frac{\lambda v_n^2}{n}Pen(f_0)<\infty$. Let $\bbeta$ be a given $p$-vector such that $\sqrt{n}\|\bbeta-\bbeta_0\|\leq c$ with $c$ a strictly positive constant. Define $\widetilde f_{n,\bbeta}=\argmax_{\{f,\,\|f\|_\infty\leq C_\infty\},\, \bbeta\in\R^p} K_{n,\lambda}(f,\bbeta).$ Then, 
\begin{eqnarray*}
v_n\|\widetilde f_{n,\bbeta}-f_0\|_n&=&\bigcirc_\P(1)\\
J(\widetilde f_{n,\bbeta}) &=&\bigcirc_\P(1).
\end{eqnarray*}
Define $\hat\bbeta_n=\argmax_{\bbeta\in\R^p} K_{n,\lambda}(\widetilde f_{n,\bbeta},\bbeta)$. Then, $$v_n\|\hat \bbeta_n-\bbeta_0\|=\bigcirc_\P(1),$$
If in addition the covariates of the linear part admit a representation of the form given in assumption $(A_{corr})$ and if the penalty satisfies (A7), then: $$\sqrt{n}\|\hat \bbeta_n-\bbeta_0\|=\bigcirc_\P(1).$$
The results still hold if the number $p$ of regression covariates goes to infinity provided the sequences $h^{1/2}p$, $n^{-\nu/(2+\nu)}p$ and $n^{-(2-\nu)/(2+\nu)}p$ go to 0 when $n$ goes to infinity.
\end{theorem}
 The proof is given in Appendix. The main keys are controls given by \cite{VanderGeer00}, relying on the entropy.

  Minimax optimality is obtained both for the linear predictor and the nonparametric estimator. This result is available for a wide class of distributions; {\it e.g.} compared to \cite{MammenVanDerGeer}, assumptions on the distributions seem weaker.  But, as discussed previously, the correlation assumption $(A_{corr})$ under which optimality is acquired, even if weaker than Rice's (1986), is more restrictive than many of those encountered in literature, and in particular comparing to \cite{MammenVanDerGeer}. Note also that without correlation conditions both estimators attain nonparametric convergence rates. 

  The fact that the results hold for a number of covariate $p$ going to infinity allows to have many covariates in the linear part. This remark can be useful for dimension reduction modeling, where the number of covariates is large. However the rate of convergence for $p$ may be poor when the regularity of the function $f$ is poor.

  In order to illustrate the general framework in which we gave the asymptotic behaviour, let us consider the penalty proposed in \cite{Antoniadis96}. To exploit the sparsity of wavelet representations, we will assume that $f$ belongs to a Besov space on the unit interval, $\mathcal B^s_{\pi,r}([0,1])$, with $s+1/\pi-1/2>0$. The last condition ensures in particular that an evaluation of $f$ at a given point makes sense. For a detailed overview of Besov spaces we refer to \cite{HardleKerkPicTsyb}.

 \begin{corollary}
Suppose $f$ belongs to a Besov ball $\mathcal B^s_{\pi,r}(C)$ with $C>0$, $s>1/2$, $0<s+1/\pi-1/2$, $\pi>2/(1+2s)$ and $1/\pi<s<\min(R,N)$, where $N$ denotes the number of vanishing moments of the wavelet $\psi$ and $R$ its regularity. 
Take the penalty: $Pen(f)=\sum_{j=j_S}^{J-1}2^{2js}\sum_k|\theta_{j,k}^W|^2$ where $\theta_{j,k}^W$ are the wavelet coefficients of $f$ and $j_S\geq 0$ a given resolution level. 
Assume conditions (A1) to (A3), (A4.1), (A7) and $(A_{corr})$ hold.\\
If $\lambda\sim n^{-2s/(1+2s)}$, we can deduce from Theorem~\ref{th1} that \begin{eqnarray*}
\|\widehat f_{n}\|_{s,2,\infty}&=&\bigcirc_\P(1),\\
n^{s/(1+2s)}\|\widehat f_{n}-f_0\|_{n}&=&\bigcirc_\P(1)\\
\sqrt{n}\|\hat \bbeta_n-\bbeta_0\|&=&\bigcirc_\P(1).
\end{eqnarray*}
where $\|f\|_{s,2,\infty}={\displaystyle \sup_{j\geq j_S}}\, 2^{j(s-1/2)}\left(\sum_k |\theta_{j,k}^W|^2\right)^{1/2}.$
\end{corollary}

\begin{proof}
\cite{BirgeMassart} establish the entropy of Besov balls $\mathcal B^s_{\pi,\infty}(1)$, with $2/(1+2s)<\pi$, is $\nu=1/s$. 
One can see that $\frac{\lambda v_n^2}{n}Pen(f)\sim\|f\|_{s,2,\infty}^2=J(f)$. Consequently the $\delta$-entropy of the functional set $\{f,J(f)\leq c\}$ can be bounded up to a constant by $\delta^{-1/s}$. We thus can apply Theorem~\ref{th1}.
\end{proof}

  One drawback of this estimation procedure is its non adaptivity: the optimal value of the smoothing parameter $\lambda$ depends on the regularity $s$ of the function, which is unknown. Actually the way the penalty term is defined by (A5) cannot lead to adaptivity since it is closely linked to the norm of the functional space where the function $f$ lies. Another penalty type may be introduced.

\subsubsection{Adaptive case}

  This section deals with the introduction of an adaptive penalty. We choose to use a $\ell^1$-penalty on the wavelet coefficients. The $\ell^1$-penalty, or LASSO, in the general least squares and likelihood settings was proposed by \cite{Tibshirani}. It is well known to lead to adaptive soft-thresholding estimators in gaussian wavelets based regressions (see \cite{DonJohn94}). We refer to \cite{AntoniadisFan} for general theory on penalization on wavelets coefficients and to \cite{LoubesVanderGeer} for the use of the $\ell^1$-penalty in functional models.

Our asymptotic results for the $\ell^1$-penalty are the following:

\begin{theorem}\label{th2}
Suppose assumptions (A1) to (A3) and assumptions (A4.2) and (A6) hold.
Suppose $f$ belongs to a Besov ball $\mathcal B^s_{\pi,r}(C)$ with $s+1/\pi-1/2>0$ and $1/2<s<\min(R,N)$, where $N$ denotes the number of vanishing moments of the wavelet $\psi$ and $R$ its regularity.

Let $\bbeta$ be a given $p$-vector such that $\sqrt{n}\|\bbeta-\bbeta_0\|\leq c$ with $c$ a strictly positive constant. 
Let $K_{n,\lambda}$ be given by equation~(\ref{critere}), with the penalty: $Pen(f)=\lambda \sum_{i=i_S}^n|\theta_i^W|$ where $(\theta_{i}^W)$ are the wavelet coefficients of $f$. Define $\widetilde f_{n,\bbeta}=\argmax_{\{f,\,\|f\|_\infty\leq C_\infty\},\, \bbeta\in\R^p} K_{n,\lambda}(f,\bbeta).$
 
Suppose $\lambda=c_0\sqrt{\log(n)}$. There exists a positive constant $C(K,\sigma_0^2)$ depending only on $K$ and $\sigma_O^2$ given in assumption (A4.2) such that if $c_0>C(K,\sigma_0^2)$, then, we have $$\left(\frac{n}{\log(n)}\right)^{s/(1+2s)}\|\widetilde f_{n,\bbeta}-f_0\|_n=\bigcirc_\P(1).$$ 
Define $\hat\bbeta_n=\argmax_{\bbeta\in\R^p} K_{n,\lambda}(\widetilde f_{n,\bbeta},\bbeta)$. Then, $$\left(\frac{n}{\log(n)}\right)^{s/(1+2s)}\|\hat \bbeta_n-\bbeta_0\|=\bigcirc_\P(1),$$
If in addition the covariates of the linear part admit a representation of the form given in assumption $(A_{corr})$, then:
$$\sqrt{n}\|\hat \bbeta_n-\bbeta_0\|=\bigcirc_\P(1).$$
The results still hold if the number $p$ of regression covariates goes to infinity provided the sequences $h^{1/2}p$, $\frac{\log(n)^{4s/(1+2s)}}{n^{(2s-1)/(1+2s)}}p$ and $\left(\frac{\log(n)}{n}\right)^{s/(1+2s)}p$ go to 0 when $n$ goes to infinity.
\end{theorem}

The proof is given in Appendix and relies on M-estimation techniques of \cite{VanderGeer00}.

  As noted before, assumption (A4.2) is more restrictive than assumption (A4.1). The results of the Theorem~\ref{th2} could probably be extended to exponential tails distributions, weakening assumption (A4.2). We refer to page 134 of \cite{VanderGeer00} (Corollaries 8.3 and 8.8) for a discussion on the price to pay to release assumption (A4.2).

  The adaptivity is acquired. As explained before, it gives the possibility of a computable procedure, without need of a cross-validation step. Yet, the parameter $\lambda$ is only chosen among an asymptotic condition and a finite sample application arises that the exact choice of this parameter is important. This will be discussed hereafter. Note also that the link with soft-thresholding is not evident, but appears in the iterative implementation of the estimators in next section.

  The optimality~of~the~functional~part~estimation~is~of~course~still~available~in~a generalized functional~model. In such models \cite{AntoniadisBesbeasSapatinas} and \cite{AntoniadisSapatinas} propose an adaptive estimation when the variance is respectively cubic or quadratic. These assumptions have to be compared with (A4.2). 
 Recently, \cite{BrownCaiZhou} introduced a method which consists of a transformation on the observations, based on the central limit theorem, in order to be able to use Gaussian framework's results. Yet, even if the asymptotic results are satisfactory, the implementation needs an important number of observations. We will see in numerical study that our procedure is easier to compute. Note also that it is available with the presence of the linear part.

\section{Algorithm and simulation study}

This section is only devoted to the adaptive $\ell^1$-type penalty. Similar algorithms are available for other penalties but a cross-validation procedure should be elaborated because of the lack of adaptivity. An advantage of the proposed estimators is that they can be easily computed, by the way of iterative algorithms. A short simulation study is also given to evaluate the performance of the estimation with finite samples.


\subsection{Algorithm}

The implementation of the estimators was performed by a {\it backfitting} algorithm, as proposed by \cite{HastieTibshirani}.

{\bf Backfitting:} Let $\bbeta^{(0)}$ be a given $p$-dimensional vector. For each iteration $k$, do:
\begin{description}
\item[Step 1:] $f^{(k+1)}=\argmax_f K_{n,\lambda}(f,\bbeta^{(k)})$.
\item[Step 2:] $\bbeta^{(k+1)}=\argmax_{\bbeta} K_{n,\lambda}(f^{(k+1)},\bbeta)$.
\end{description}
The algorithm is stopped either when a maximal number of iterations $\kappa$ is attained or when the algorithm is stabilized, {\it i.e.} when $\|\bbeta^{(k)}-\bbeta^{(k-1)}\|\leq \delta\|\bbeta^{(k-1)}\|$ for a given tolerance value $\delta$. The returned values are $\hat\bbeta_n=\bbeta^{(K)}$ and $\hat f_n=f^{(K)}$ with $K$ maximal number of iterations of the algorithm.

To compute each of the two steps we apply a classical Fisher-scoring algorithm, detailed among others page 40 of \cite{MacCullaghNelder}. Usual in generalized models, this algorithm consists of building new variables of interest by a gradient descending method, in order to apply a ponderate regression on these new variables. Recall the notations of the GPLM given in equation~(\ref{GPLM}): one has $\eta(X,t)=X^T\bbeta+f(t)$ and $\mu(X,t)=\dot b(\eta(X,t))$. We will omit the dependence to $(X,t)$ for the sake of simplicity.

\begin{description}
\item[{\bf Step 1}:] ~\\ Note $f^{(k,0)}=f^{(k)}$.
Repeat the following iteration for $j=0\dots J_1-1$:\\
\begin{itemize}
\item[] Let $\eta^{(k,j)}=X\bbeta^{(k)}+f^{(k,j)}$. 
\item[] We introduce $Y^{(k,j)}=f^{(k,j)}+(y-\mu^{(k,j)})\left.\frac{d \eta}{d \mu}\right|_{\mu=\mu^{(k,j)}}$ and $W^{(k,j)}=\text{diag}\left(\left.\frac{d\eta}{d\mu}\right|_{\mu=\mu^{(k,j)}}\right)$. 
\item[] With an $\ell^1$-penalty, we establish that $f^{(k,j+1)}$ is a nonlinear wavelet estimator for the observations $Y^{(k,j)}$, obtained by soft-thresholding of the wavelet coefficients (see \cite{DonJohn94,DonJohn95,DonJohn98}). The threshold levels are $\lambda\Psi {W^{(k,j)-1}}\Psi^T\1_{n\times 1}$ where $\Psi$ denotes the forward wavelet transform and $\Psi^T$ the inverse wavelet transform.
\end{itemize}
Take $f^{(k+1)}=f^{(k,J_1)}$.\\

\item[{\bf Step 2:}]~\\
 Define $\bbeta^{(k,0)}=\bbeta^{(k)}$. Repeat the following iteration for $j=0\dots J_2-1$:\\
\begin{itemize}
\item[] Let $\widetilde \eta^{(k,j)}=X\bbeta^{(k,j)}+f^{(k+1)}$. 
\item[] We introduce $\widetilde Y^{(k,j)}=X\bbeta^{(k)}+(y-\widetilde \mu^{(k,j)})\left.\frac{d \eta}{d \mu}\right|_{\mu=\widetilde \mu^{(k,j)}}$ and $\widetilde W^{(k,j)}=\left(\left.\frac{d\eta}{d\mu}\right|_{\mu=\widetilde \mu^{(k,j)}}\right)$. 
\item[] Then $\bbeta^{(k,j+1)}$ is the regression parameter of $Y^{(k,j)}$ on $X$ with ponderations $\widetilde W^{(k,j)}$, {\it i.e.} $\bbeta^{(k,j+1)}=(X^T  {\widetilde W}^{(k,j){-1}}X)^{-1}X^T {\widetilde W}^{(k,j){-1}}\widetilde Y^{(k,j)}$.
\end{itemize}
Take $\bbeta^{(k+1)}=\bbeta^{(k,J_2)}$.
\end{description}
Maximal number of iterations $J_1$ and $J_2$ can be fixed to $1$ to simplify the algorithm (this is what is proposed actually in \cite{Muller01}). In computation studies, no main difference is observed while modifying parameters $J_1$ and $J_2$. We therefore also decide to fix them equal to 1. Note that the only matrix needing inversion is diagonal so we can expect a fast computation. 

The initialization values are set as follows: for all~$j=1,\ldots,p$, $\bbeta^{(0)}_j=0$ and $f\equiv 0$, except for the Poisson distribution were for all $i=1,\ldots,n$, $f^{(0)}(t_i)=G(y_i)$, with $G$ the link function associated to the model (with a slight modification if the value does not exist).

In the particular Gaussian case, the two steps are non iterative. We explicit how the estimators are implemented in a Gaussian framework to ameliorate the comprehension of the algorithm:
\begin{description}
\item[{\bf Step 1:}] The iterate $f^{(k+1)}$ is the wavelet estimator for the observations $y-X^T\bbeta^{(k)}$, with a soft-thresholding on wavelet coefficients with an uniform threshold level $\lambda$.
\item[{\bf Step 2:}] We obtain $\bbeta^{(k+1)}$ by a maximum likelihood estimation on  $y-f^{(k+1)}$; this means that one has $\bbeta^{(k+1)}=(X^T X)^{-1}X^T (Y-f^{(k+1)})$.
\end{description}
We can recognize the {\it backfitting} algorithm studied in \cite{ChangQu}, \cite{FadiliBullmore} and \cite{gannaz}. In a Gaussian framework, the variance in observations is constant and consequently the threshold level is uniform. In a generalized framework, the matrix $W$ ponderates the threshold level to take into account the inhomogeneity of the variance.

In generalized functional models, {\it i.e.}~without the presence of a linear part, the numerical implementation of the estimators proposed here has already been explored by \cite{Sardy}. The authors propose an interior point algorithm based on the dual maximisation problem. Comparing to this resolution scheme, our procedure has the advantage of an easy interpretation of the different steps in the algorithm. As noted previously, wavelet estimators have also been explored by \cite{AntoniadisBesbeasSapatinas}, \cite{AntoniadisSapatinas} and \cite{BrownCaiZhou}. These papers need to aggregate the data into a given number of bins. If the two first cited papers allow to choose small size of bins, it appears to be quite a constraint for the third one. Actually, it is worthy noticing the simplicity of our algorithm, {\it e.g.}~comparing with those implementations.

Finally we can remark that the algorithm establishes the link between the soft-thresholding procedure and the $\ell^1$-penalty, usual in Gaussian models. We easily can see that other penalties will lead to the usual thresholding scheme they are associated with. Following \cite{AntoniadisFan} we can extend this algorithm introducing other thresholding schemes {\it e.g.} to hard-thresholding which will correspond to an $\ell^0$ penalty, or to SCAD-thresholding which will be given by a mixed $\ell^0$ and $\ell^1$ penalty (see \cite{FanLi}).

\subsection{Simulations}

In this subsection, we give some simulation results. All the calculations were carried out in R version 2.14.1 \nocite{R} on a Unix environment. For the DWT, we used the {\it wavelets} package 0.2-6 developed by \cite{Rwavelets}. In order to better study the quality of our estimates, we compared with kernel and splines based estimators. The computation of those procedures was done using the {\it KernGPLM} package developped by \cite{RKernGPLM} and applied {e.g.} in \cite{Muller01}. This package proposes four estimation procedures, based on kernel or splines estimators and using backfitting or Speckman's type algorithms. They are available for Gaussian and Bernoulli distributions and thus they have been extended here to Binomial and Poisson frameworks.

A cross-validation procedure has been implemented for kernel and splines methods. It is based on the Generalized Cross Validation criterion initially proposed by \cite{CravenWahba}. It has also been applied by \cite{SullivanYandellRaynor} in generalized functional regression models and suggested by \cite{Speckman} for partially linear models. Let $\hat \mu_n=G^{-1}(\bX^T\hat\bbeta_n+\hat f_n)$ be the resulting estimation of the means and $\hat v_n=\phi \ddot b(\bX^T\hat\bbeta_n+\hat f_n)$ the estimation of the variances. If $\hat \mu_n=\mathcal R(\lambda)y$, with $\mathcal R(\lambda)$ operator depending on the degree of smoothness $\lambda$ of the procedure, than the GCV score is equal to $$GCV(\lambda)=\frac{\sum_{i=1}^n (y_i-\hat\mu_i)^2/\hat v_i}{\left(n-trace(\mathcal R(\lambda))\right)^2}.$$
The denominator $n-trace(\mathcal R(\lambda))$ corresponds to the number of degree of freedom of the model. It is used for testing the GPLM by \cite{HardleMammenMuller} or \cite{Muller01}. Its computation is detailed in \cite{Muller01} and is available in the R-package of \cite{RKernGPLM}. In this paper we restrict the GCV procedure to kernel and splines based estimators but the GCV can also be proposed for an automatic choice of the wavelet thresholds (see {\it e.g.} \cite{JansenMalfaitBultheel}).  Note that \cite{FadiliBullmore} used the GCV procedure combined with the $\ell^2$ wavelets penalty of Corallary 1 in Gaussian partially linear models. This estimation scheme could be extended here to others distributions using the algorithm of previous section. 

We simulate $n = 2^8$ observations. The covariates are written according to assumption $(A_{corr})$: $X_i=(x_{1,i},\,x_{2,i})$ with \begin{itemize}
\item[$\bullet$] $x_{1,i}=g_1(i/n)+\xi_{1,i},$ with $\xi_{1,i}$ independent and identically distributed variables following a uniform distribution on the interval $[-0.5,0.5]$ and with $g_1(x)=2 x-1$. 
\item[$\bullet$] $x_{2,i}=g_2(i/n)+0.5\xi_{2,i},$ with $\xi_{2,i}$ independent and identically distributed variables following a standard normal distribution and with $g_2(x)=(x-0.5)^2$.
\end{itemize}
The parameter $\bbeta$ was taken equal to $(1,\,-1)$.

\begin{figure}[!h]
  \includegraphics[scale=0.8]{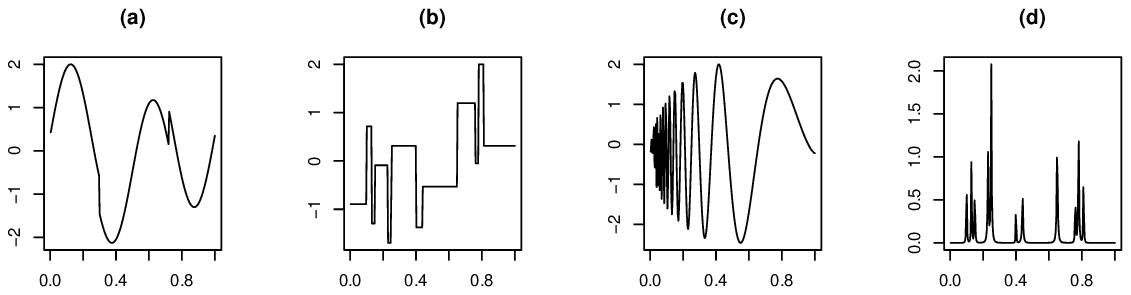}
\caption{Functional part of the generalized partially linear model. Figure~(a) corresponds to the {\it Heavisine} function, Figure~(b) to the {\it Blocks} function, Figure~(c) to the {\it Doppler} function and Figure~(d) to the {\it Bumps function}.}
\label{fig:functions}       
\end{figure}

The four test signals introduced by \cite{DonJohn94} were considered in the functional part~$f$: \begin{itemize}\item the {\it Heavisine} function, for a smooth function, combination of sinusoidal functions, \item the {\it Blocks} function, for a piecewise constant function, \item the {\it Doppler} function, for a high frequency signal \item and the {\it Bumps} function for a function presenting localised high variations.\end{itemize} These functions are given in Figure~\ref{fig:functions}. 

We will more precisely study the estimation quality for 
\begin{itemize}
\item a Gaussian distribution; observations are $y_i\sim\mathcal N(\eta_i, \sigma^2)$, with $\sigma^2=\phi$,
\item a Binomial distribution; observations are $y_i$ such that $y_i\times m\sim\mathcal B(\mu_i, m)$ with $m=\phi^{-1}$ and the logit link function, {\it i.e.} $\mu_i=\frac{\exp(\eta_i)}{1+\exp(\eta_i)}$,
\item a Poisson distribution; observations are $y_i\sim\mathcal P(\mu_i)$, with the link function, {\it i.e.} $\mu_i=\exp(\eta_i)$,
\end{itemize} 
as these distributions seem to be the most frequently encountered in modelization.

To conclude with respect to the quality of estimation, we will give the mean value of the squared error for the parameter $\bbeta$, noted $MSE_{\bbeta}$ which is the empirical mean of $SE_\beta=\sum_{j=1}^p (\hat\bbeta_j-\bbeta_j)^2$. For the nonparametric part, we will evaluate the average mean squared error, noted $AMSE_f$, which is the empirical mean of $MSE_f=\frac 1 n \sum_{i=1}^n\left(\hat f_n(t_i)-f_0(t_i)\right)^2$. To evaluate the global estimation quality, we will also compute a global average mean sqared error $AMSE$, empirical mean of $MSE=\sum(\frac 1 n \sum_{i=1}^n\left( \bX_i^T\hat\bbeta_n+\hat f_n(t_i)-\bX_i^T\bbeta_0-f_0(t_i)\right)^2.$ 

All results were obtained on $500$ simulations with the same covariates $X_i$ and the same functional parts~$f$. A maximal number of $\kappa=200$ iterations was taken and the tolerance value defined above was equal to $\delta=10^{-10}$ when applying the algorithms. The Daubechies's wavelets base with a filter length of 8 was chosen. Concerning the kernel estimators, the Biweight kernel was used, with a bandwidth varying from 0.005 to 0.05 with a step of 0.005. The splines estimators were computed with a smoothing parameter varying from 0.2 to 0.7 with a step of 0.05. Kernel bandwidths and splines smoothing parameters are chosen by minimizing the GCV score.

\subsubsection{Preliminary study: choice of the threshold level}

\begin{figure}[!h]
 \centering{\includegraphics[width=150mm]{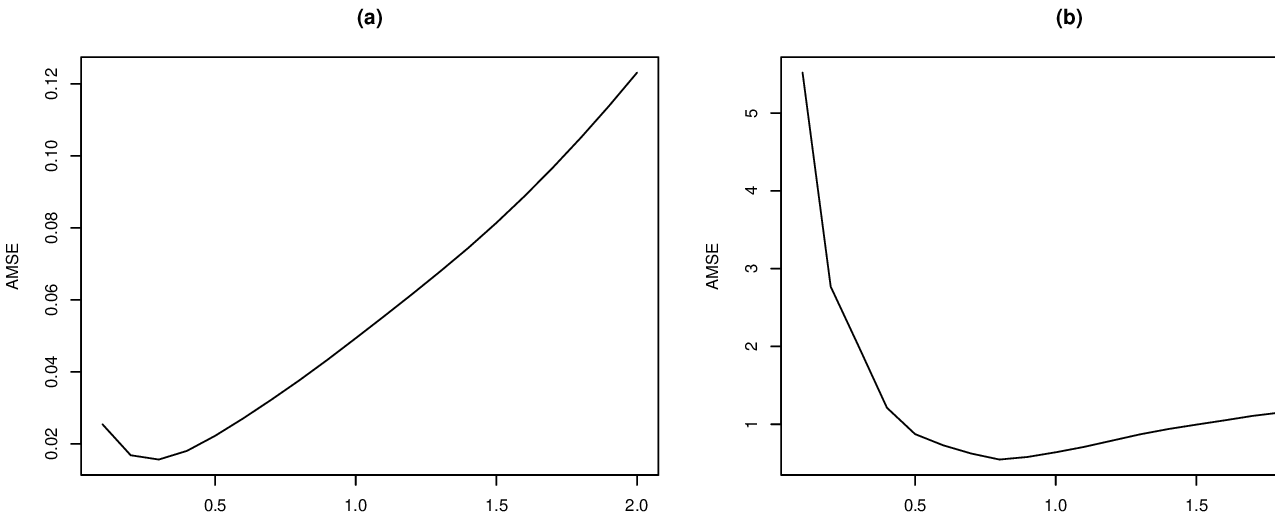}}
\caption{Evolution of the functional AMSE with respect to the threshold level in a generalized functional model with the {\it Heavisine} function. Figure (a) corresponds to a Gaussian distribution with the dispersion parameter $\phi=0.05$, Figure (b) to a Binomial distribution with the dispersion parameter $\phi=0.05$ and Figure (c) to a Poisson distribution. Calculations were done on 50 data sets of size $n=2^8$. }
\label{fig:evolution}       
\end{figure}

The asymptotic behaviour of the estimators only defines the threshold level $\lambda$ up to a constant. Yet, numerical implementation needs to determine the exact threshold level $\lambda$ in the algorithm. In practice, one can see it has an important impact on the quality of estimation. Figure~\ref{fig:evolution} gives the evolution of the AMSE with respect to the threshold level in the GPLM with different distributions. One can observe that the threshold levels attaining the minima are very different among the distribution.

With a Gaussian distribution,~following \cite{DonJohnKerkPic}, we choose $\lambda=\sqrt{2\phi\log(n)}$. This choice overevaluates the optimal threshold level in many cases but is the most often encountered in practice.


\paragraph{Binomial distribution.}


Due to homogeneity reasons, it seems well-adapted to fix a threshold level of the form $\lambda=c'_0\sqrt{\phi\log(n)}$, where $\phi$ corresponds to the dispersion parameter in distribution (\ref{modele}) and $c'_0$ denotes a positive constant.

\begin{figure}[!h]
 \begin{center} \includegraphics[width=120mm,height=50mm]{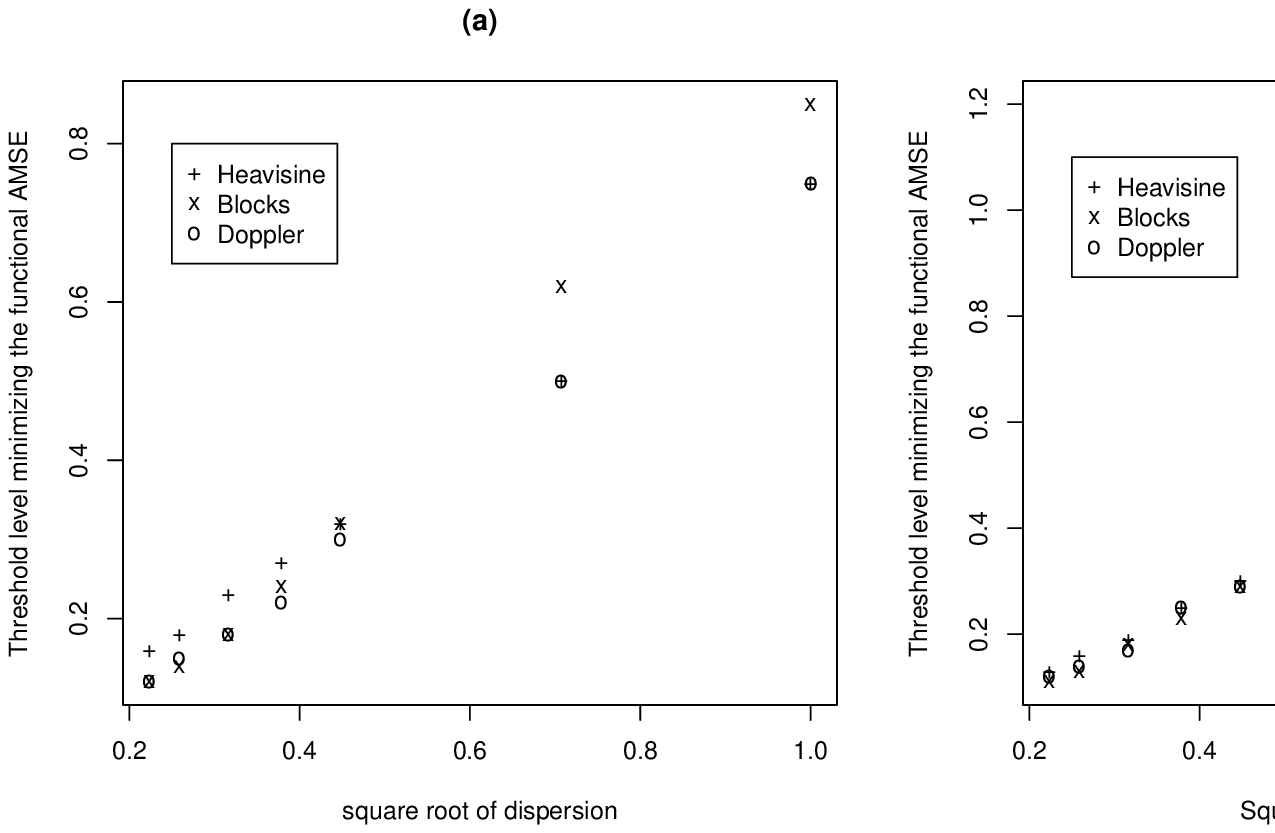}
\end{center}
\caption{Evolution of the threshold minimizing the functional AMSE when estimating the {\it Heavisine}, {\it Blocks} and {\it Doppler} functions in a Binomial setup, with respect to $\sqrt{\phi}$. In a GFM on Figure~(a) and in a GPLM in Figure~(b). Calculations were done on $100$ samples of size $n=2^8$.}
\label{fig:seuils_bino}       
\end{figure}

To ensure this conjecture in a Binomial setup, we plot the optimal threshold obtained for different value of the Binomial parameter $m$. The study was done on a generalized functional model (GFM) and on a GPLM. Figure~\ref{fig:seuils_bino} confirms that the chosen form seems appropriate. Linear fittings are given in Table~\ref{tab:seuils_binomial}. The {\it Bumps} signal was not considered during this preliminary step because the evaluation of the minima in the functional AMSE leads to a threshold level which clearly over-smoothes the function.

\vspace{0.2cm}

\begin{table}[h]

\caption{Numerical indexes for the regression of the threshold level that minimizes the functional $AMSE$ with respect to $\sqrt{\phi \log(n)}$ with a Binomial distribution. Calculations were done on $100$ samples of size $n=2^8$.}

\label{tab:seuils_binomial}
\begin{center}
\begin{tabular}{cc}

\hspace{0.5cm}

\begin{minipage}{0.5\textwidth}{\begin{tabular}{@{} l@{\hspace{0.5cm}} lll @{}}
\hline\noalign{\smallskip}
  \multicolumn{4}{c}{Generalized functional model}\\
\noalign{\smallskip}\hline\noalign{\smallskip}
  Function & {\it Heavisine} & {\it Blocks} & {\it Doppler}\\
  $R^2$ coefficient & 0.999 & 0.992 & 0.988\\
 constant $c'_0$ & 0.306 & 0.269 &  0.258\\
\noalign{\smallskip}\hline
\end{tabular}}\end{minipage}

& 

\begin{minipage}{0.5\textwidth}{\begin{tabular}{@{} l@{\hspace{0.5cm}} lll @{}}
\hline\noalign{\smallskip}
  \multicolumn{4}{c}{Generalized partially linear model}\\
\noalign{\smallskip}\hline\noalign{\smallskip}
  Function & {\it Heavisine} & {\it Blocks} & {\it Doppler}\\
  $R^2$ coefficient & 0.997 & 0.999 & 0.988\\
 constant $c'_0$ & 0.273 &  0.257 & 0.259\\
\noalign{\smallskip}\hline
\end{tabular}}
\end{minipage}

\end{tabular}
\end{center}
\end{table}

The linear fitting is coherent provided the $R^2$ coefficients. The constant $c'_0$ obtained by linear fitting is varying with respect to the simulated samples, but the observed variations are small. The mean value is $0.306$ and the standard deviation is $0.019$. The conjecture of a uniform constant seems acceptable. We therefore choose to take the value $c'_0=0.3$ in the following for Binomial distributions. Note that due to the central limit theorem, one would have expected to take $c'_0=\sqrt{2}$ for large values of the parameter $m$, but our simulation study shows that this would lead to an oversmoothing for small values of $m$.

Note all the calculations were done with a fixed sample size $n=2^8$. To better evaluate the form of the optimal threshold in practice, one could also study the evolution of the threshold level with respect to the sample size $n$.


\paragraph{Poisson distribution.} Estimation in a Poisson functional model has been more intensively explored. Note that \cite{Sardy} propose a threshold level. The main drawback is that the level given depends on the estimated function. Yet, the choice is based on an universal large deviation inequality which does not seem to be well-adapted in this procedure. Indeed, due to the iterative interpretation of the estimation, the inhomogeneity of the variance of the observations is taken into account within the estimation. 

Recently, \cite{ReynaudRivoirard} have developed a procedure based on wavelet hard-thresholding estimation for Poisson regression. In their estimation the thresholding step is defined directly and not through a penalization procedure like here. The authors then present a detailed numerical study showing the high instability of the optimal threshold level with respect to the estimated function. 

\begin{figure}[!h]
\centering{    \includegraphics[width=60mm]{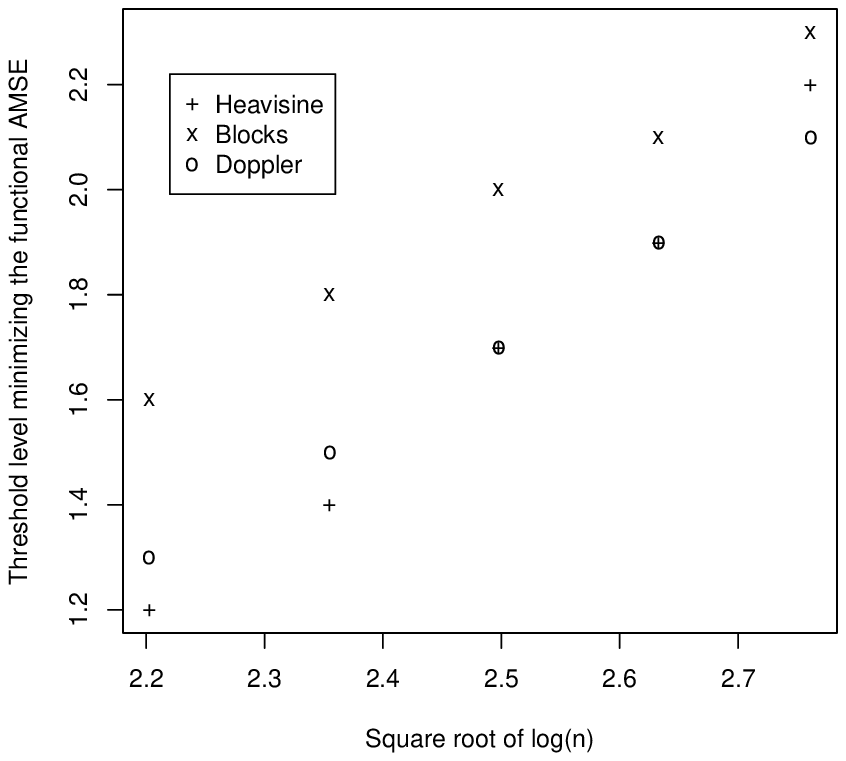}}
\caption{Evolution of the optimal threshold with respect to $\sqrt{\log(n)}$ with a Poisson distribution in a Generalized functional model. Calculations were done on $100$ samples for each value of $n$.}
\label{fig:seuils_poisson}       
\end{figure}

In our procedure, we observe that the optimal threshold level does not vary significantly with respect to the estimated function; this can bee seen on Figure~\ref{fig:seuils_poisson} which represents the evolution of the threshold level which minimizes the AMSE with respect to $\sqrt{\log(n)}$ in a Poisson functional regression. The minima can clearly be identified in the generalized functional model but they are much more approximative in a partially linear context. Consequently, we prefer not to take them into account in our study.


\begin{table}[h]

\caption{Numerical indexes for the regression of the threshold level that minimizes the functional $AMSE$ with respect to $\sqrt{\log(n)}$ with a Poisson distribution. Calculations were done on $100$ samples for each size of $n$.}

\label{tab:seuils_poisson}
\begin{center}

\hspace{0.5cm}

\begin{tabular}{@{} l@{\hspace{0.5cm}} lll @{}}
\hline\noalign{\smallskip}
  \multicolumn{4}{c}{Generalized functional model}\\
\noalign{\smallskip}\hline\noalign{\smallskip}
  Function & {\it Heavisine} & {\it Blocks} & {\it Doppler}\\
  $R^2$ coefficient & 0.996 & 0.990 & 0.989\\
 constant $c'_0$ & 0.682 & 0.804 &  0.679\\
\noalign{\smallskip}\hline
\end{tabular}



\end{center}
\end{table}


As the theoretical result is asymptotic, we should study larger values of the sample size $n$. According to the results in Table~\ref{tab:seuils_poisson}, one may choose a constant approximatively equal to $0.72$ in a functional model and in a GPLM. We therefore choose to take the value $c'_0=0.72$ in the following for Poisson distribution.

\newpage

\paragraph{Remark:} In a Gaussian or a Binomial regression, one may need the dispersion parameter~$\phi$. Actually in literature, it is classically estimated at each iteration by $$\phi^{(k)}=\frac{1}{n}\sum_{i=1}^n\frac{(y_i-\mu_i^{(k)})^2}{\ddot b(\eta^{(k)})}.$$ Due to the bad quality of this estimator in GPLM, we prefer to consider in this paper that the dispersion parameter is known. In a Gaussian model, \cite{gannaz} proposed an efficient QR-based estimator for $\phi$. It would be interesting to explore whether it could be extended to generalized models.

\subsubsection{Example 1: Gaussian distribution}

Example 1 deals with a Gaussian model. The Gaussian distribution implementation is not a novelty for this estimation procedure, and we refer to \cite{ChangQu}, \cite{FadiliBullmore} and \cite{gannaz} for detailed studies on simulated or real values data. We briefly consider this case in order to have a comparison base for other distributions. 

The signal-to-noise ratio (SNR) of a signal is defined as the norm of the ratio of the mean value with respect to the standard deviation. In GPLM the SNR for the nonparametric part, noted $SNR_f$ and the SNR for the linear part of the model, noted $SNR_{\bbeta}$, are respectively equal to \begin{eqnarray*} SNR_f^2&=&\frac{1}{n}\sum_{i=1}^n \frac{f_0(t_i)^2}{\phi\,\ddot b\left(\bX_i^T\bbeta_0+f_0(t_i)\right)},\\
\text{and~~~}SNR_{\bbeta}^2&=&\frac{1}{n}\sum_{i=1}^n \frac{(\bX_i^T\bbeta_0)^2}{\phi\,\ddot b\left(\bX_i^T\bbeta_0+f_0(t_i)\right)}.
\end{eqnarray*} With a high SNR, say approximatively $5$, one can expect a good quality of estimation, while with a small value, like $1$, the quality of estimation cannot be satisfying.

\begin{figure}[!h]
\centering{   \includegraphics[width=150mm,height=45mm]{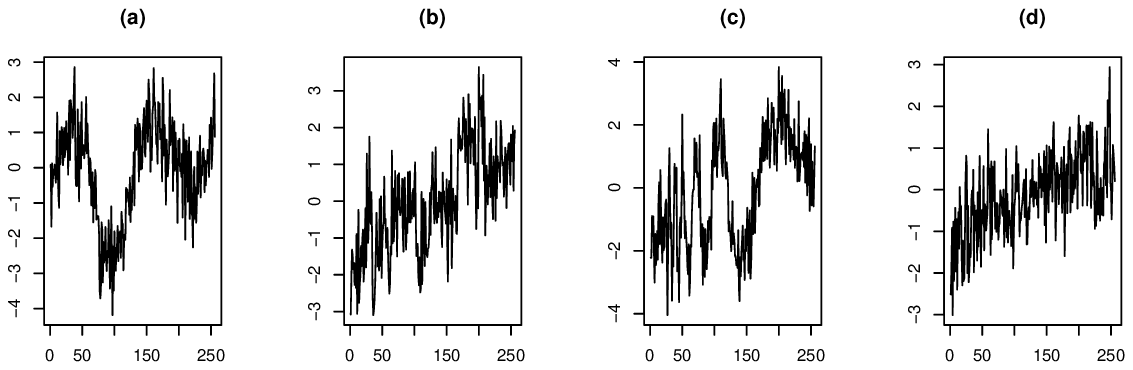}}
\caption{An example of a simulated data set in Example 1, with the {\it Heavisine} function in Figure (a), the {\it Blocks} function in Figure (b), the {\it Doppler} function in Figure (c) and the {\it Bumps} function in Figure (d).}
\label{fig:obs_gplm_gauss}       
\end{figure}

Example 1 considered Gaussian observations with a variance $\phi=\sigma^2=0.05$. An example of simulated observations obtained in this example are given in Figure~\ref{fig:obs_gplm_gauss}. The $SNR_{\bbeta}$ for the linear regressor is approximatively equal to $3.8$ with each target signals and the functional $SNR_f$ is respectively equal to $5.5$ for the {\it Heavisine} function, $3.6$ for the {\it Blocks} function, $5.8$ for the {\it Doppler} function and $0.97$ for the {\it Bumps} function. Except for the $SNR_f$ of the {\it Bumps} signal, the SNRs are high and we expect a good quality in estimation. 

In order to evaluate if the quality observed for the wavelets based procedure is due or not to the presence of the linear part, we give in Table~{\ref{tab:table_gauss_nonpara} the $AMSE$ for an usual functional regression model $y_i\sim\mathcal N\left(f(t_i),0.05\right)$, {\it i.e.} without linear part.

\begin{table}[h]

\caption{Measures of quality the wavelet estimator over the 500 simulations in an gaussian functional regression model $y_i=f(t_i)+\epsilon_i$ with $n=2^{8}$, $\epsilon_i$ following a gaussian distribution $\mathcal N(0, 0.05)$ and differents functions $f$.}

\label{tab:table_gauss_nonpara}

\begin{center}

{\begin{tabular}{@{} l@{\hspace{1cm}} l @{\hspace{1cm}} l @{\hspace{1cm}} l @{\hspace{1cm}} l @{}}
\hline\noalign{\smallskip}
  $f$  & Heavisine & Blocks & Doppler & Bumps\\
\noalign{\smallskip}\hline\noalign{\smallskip}
 $AMSE_f$ & 0.03411 & 0.01061 & 0.698 & 0.0403\\
  
 \noalign{\smallskip}\hline
\end{tabular}}

\end{center}

\vspace{-12pt}

\end{table}

All the numerical measures for each estimation procedure are given in Table~\ref{tab:table_gauss}. This Table is completed by the boxplots of the squared errors for the linear part MSE, the mean squared errors for the functionnal part and the global mean squared errors, in Figure~\ref{fig:boxplot_gauss}.

The quality indexes show that the wavelets based estimators do not perform as well as the kernel or spline based methods when estimating the functions {\it Heavisine} or {\it Blocks}. The quality is also lower for the linear part when considering such functional parts. Yet, we can see with Table~\ref{tab:table_gauss_nonpara} that even in a functional model the AMSE for wavelets estimators are lower than those obtained by kernel or splines estimators. It would be interesting in such contexts to study if a GCV step improves the quality of the wavelets-based procedure.

Concerning the simulations with the {\it Doppler} or the {\it Bumps} functions, we can see that wavelets-based procedure is satisfactory. The best AMSE for the functional part is given by the kernel method with Speckman algorithm but then wavelets give results comparable with the others methods. Note also that the boxplots illustrate that wavelets lead to a more stable estimation of the linear part.

We can remark that in this example, whatever the functional part is, the backfitting procedures perform worse than the Speckman algorithms. The kernel based estimator with Speckman algorithm gives in many cases the best quality indexes. If results are close to the ones of splines estimator with Speckman algorithm for {\it Heavisine} and {\it Blocks} signals, it is more performant for {\it Doppler} and {\it Bumps} nonparametric parts.

Finally, to illustrate the visual quality of our functional estimators, we give in Figure~\ref{fig:gplm_gauss} the estimated functions for one simulation with each test functions. As expected, the estimates tends to oversmooth the functions (see \cite{DonJohn94}). Other threshold levels were proposed in literature but not tested here. The peaks in the {\it Bumps} signal are not all identified by our procedure but this is coherent with the small SNR.

\vspace{-12pt}

\begin{figure}[!h]
\centering{   \includegraphics[width=150mm,height=45mm]{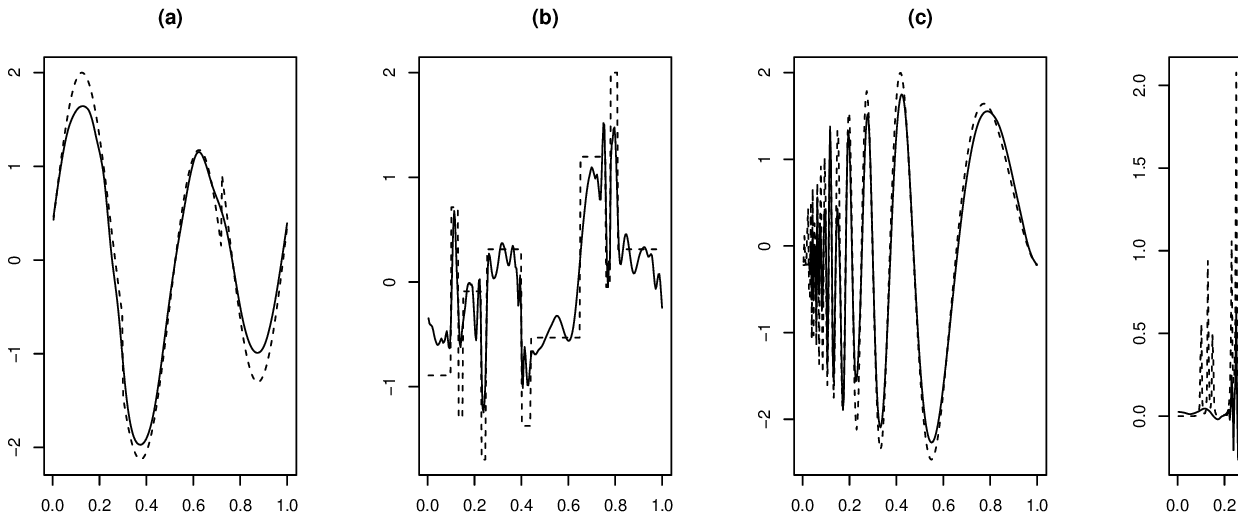}}
\caption{An example of estimation of the nonparametric part in Example 1, with the {\it Heavisine} function in Figure (a), the {\it Blocks} function in Figure (b), the {\it Doppler} function in Figure (c) and the {\it Bumps} function in Figure (d). Dots lines corresponds to the true functions and plain lines to their estimate. }
\label{fig:gplm_gauss}       
\end{figure}

\begin{table}[!h]

\caption{Measures of quality the estimates over the 500 simulations in Example 1 with $n=2^{8}$ and $\phi=0.05$. Lowest values are in bold face type.}

\label{tab:table_gauss}

\begin{center}
\vspace{-12pt}

\begin{tabular}{cc}
\begin{minipage}{0.5\textwidth}{
\begin{tabular}{@{} l@{\hspace{0.5cm}} l @{\hspace{0.5cm}} l @{\hspace{0.5cm}}l @{}}

\hline\noalign{\smallskip}

&\multicolumn{3}{c}{{\it Heavisine}} \\ [5pt]
   & {$MSE_{\bbeta}$} & {$AMSE_f$} & Global AMSE\\
\noalign{\smallskip}\hline\noalign{\smallskip}
  {\it Wavelets} &  0.0463 & 0.0531&0.0342\\
  {\it Kernel Speckman} & {\bf 0.0033} & 0.0104&0.0099\\
  {\it Kernel Backfitting} & 0.0082 & 0.0164& 0.0150\\
 {\it Spline Speckman} & 0.0040 & {\bf 0.0104} & {\bf 0.0095}\\
  {\it Spline Backfitting} & 0.0080 & 0.0157 & 0.0143\\ 
 
 \noalign{\smallskip}\hline \\ 

&\multicolumn{3}{c}{{\it Blocks}} \\ [5pt]
   & {$MSE_{\bbeta}$} & {$AMSE_f$} & Global AMSE\\
\noalign{\smallskip}\hline\noalign{\smallskip}
  {\it Wavelets} & 0.0616 & 0.1465 & 0.1217\\
  {\it Kernel Speckman} & {\bf 0.0064} & {\bf 0.0398} &{\bf 0.0382}\\
  {\it Kernel Backfitting} & 0.0109 & 0.0577 &0.0550\\
 {\it Spline Speckman} & 0.0069 & 0.0438 & 0.0423\\
  {\it Spline Backfitting} & 0.0152 & 0.0603 & 0.0561\\ 

 \noalign{\smallskip}\hline

\end{tabular}}\end{minipage}
 &
\begin{minipage}{0.5\textwidth}{
\begin{tabular}{@{} l@{\hspace{0.5cm}} l @{\hspace{0.5cm}} l @{\hspace{0.5cm}}l @{}}

 \hline \noalign{\smallskip}

&\multicolumn{3}{c}{{\it Doppler}} \\ [5pt]
   & {$MSE_{\bbeta}$} & {$AMSE_f$} & Global AMSE\\
\noalign{\smallskip}\hline\noalign{\smallskip}
  {\it Wavelets} &  0.0073 & 0.0700 & 0.0685\\
  {\it Kernel Speckman} & 0.0094 & {\bf 0.0473} & {\bf 0.0450}\\
  {\it Kernel Backfitting} & 0.0093 & 0.0721 & 0.0700\\
 {\it Spline Speckman} & {\bf 0.0071} & 0.0757 & 0.0737\\
  {\it Spline Backfitting} & 0.0125 & 0.0920 & 0.0884\\

 \noalign{\smallskip}\hline \\

&\multicolumn{3}{c}{{\it Bumps}} \\ [5pt]
   & {$MSE_{\bbeta}$} & {$AMSE_f$} & Global AMSE\\
\noalign{\smallskip}\hline\noalign{\smallskip}
  {\it Wavelets} & {\bf 0.0015} & 0.0378 & 0.0381\\
  {\it Kernel Speckman} & 0.0051 & {\bf 0.0318} &{\bf 0.0380}\\
  {\it Kernel Backfitting} &0.0092 & 0.0393 &0.0375\\
 {\it Spline Speckman} & 0.0040 & 0.0334 &0.0329\\
  {\it Spline Backfitting} & 0.0071 & 0.0398 &0.0388\\
 
\noalign{\smallskip}\hline 
 
\end{tabular}}
\end{minipage}

\end{tabular}

\vspace{12pt}

\end{center}

\end{table}


\begin{figure}[h]

\centering{
{\it Squared Error for $\bbeta$}\\
  \includegraphics[width=150mm,height=35mm]{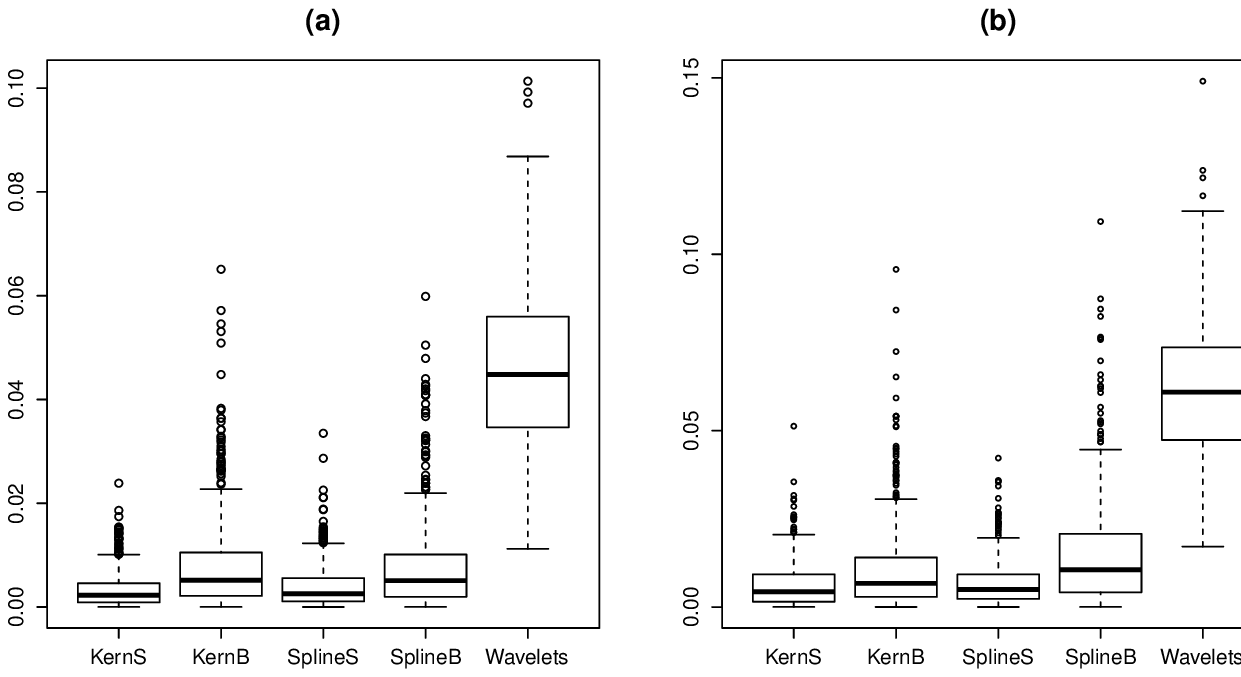}\\
{\it Mean Squared Error for the functional part }\\ 
 \includegraphics[width=150mm,height=35mm]{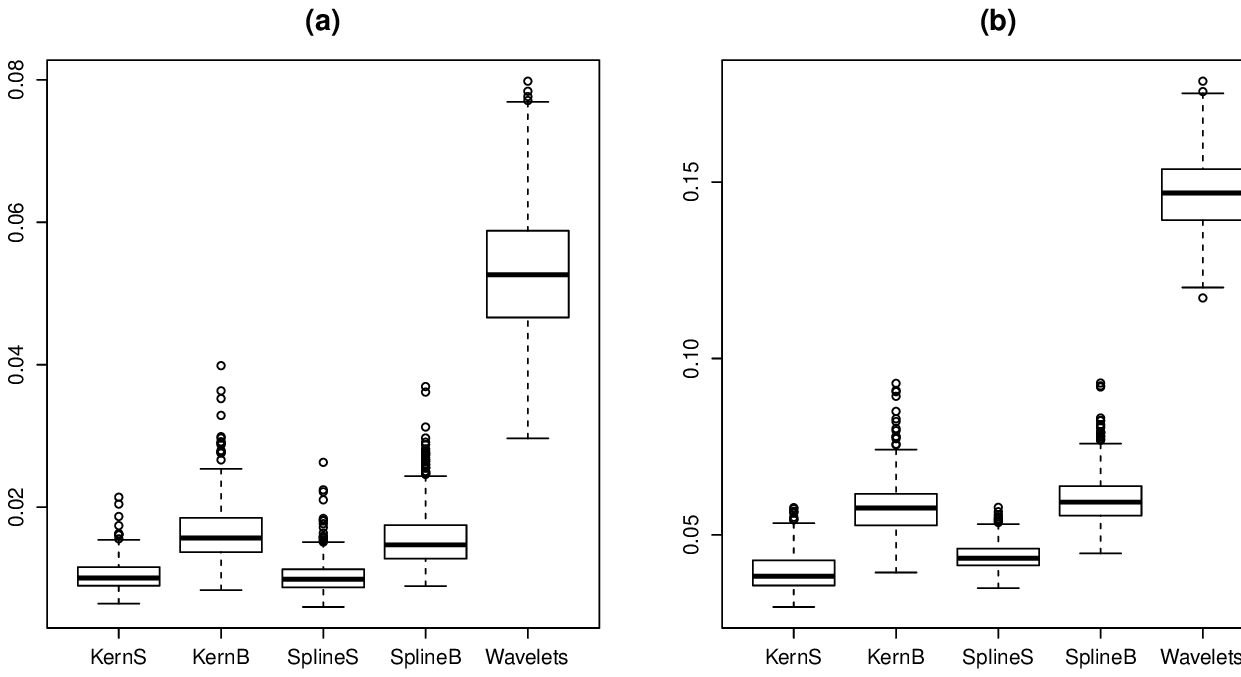}\\
{\it Global Mean Squared Error}\\
 \includegraphics[width=150mm,height=35mm]{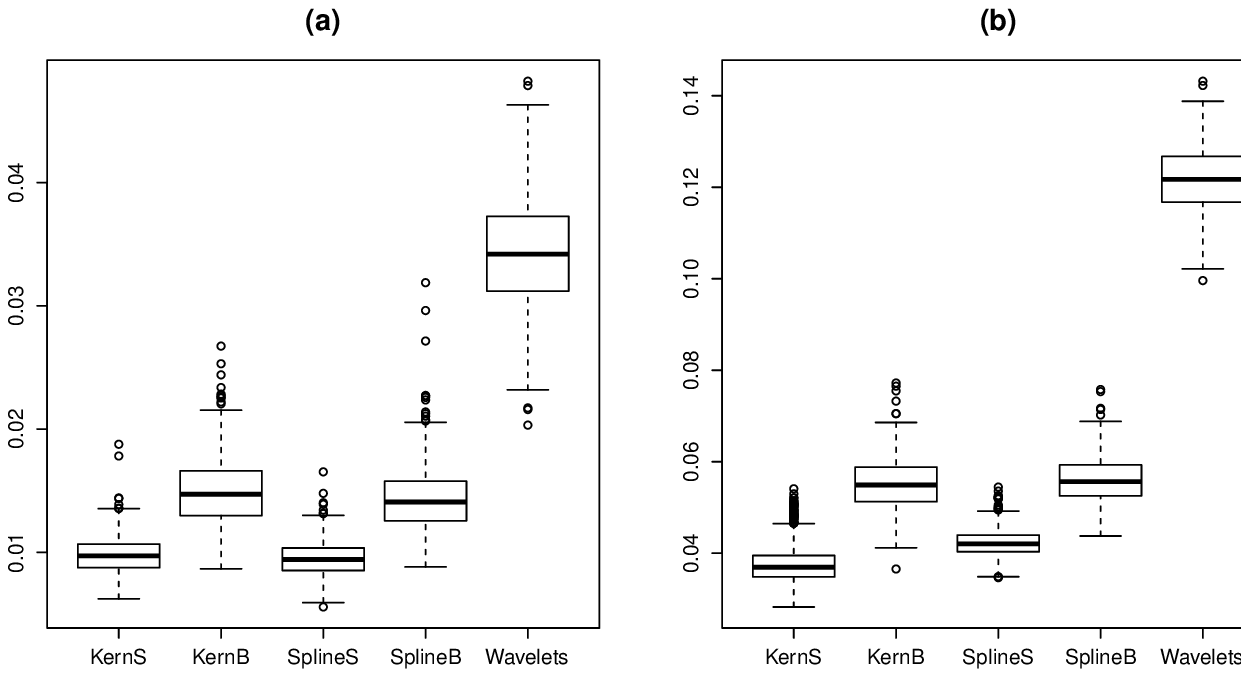}\\
}
\caption{Boxplots of the MSE in Example 1, with the {\it Heavisine} function in Figures (a), the {\it Blocks} function in Figures (b), the {\it Doppler} function in Figures (c) and the {\it Bumps} function in Figures (d). {\it KernS} and {\it KernB} stand for the kernel procedures respectively with Speckman algorithm and Backfitting algorithm, and {\it SplineS} and {\it SplineB} stands for the splines procedures respectively with Speckman algorithm and Backfitting algorithm.}
\label{fig:boxplot_gauss}   
    
   \vspace{-3cm}
    
\end{figure}

\newpage

~\\
\newpage



\subsubsection{Example 2: Binomial distribution}

In Example 2, we consider a Binomial distribution: observations $Y_i$ are such that $Y_i\times m$ are independently drawn from Binomial distributions $\mathcal B(\mu_i, m)$ with the parameter $m$ equal to $m=\phi^{-1}$. The link function considered is the logistic. The mean is thus equal to $$\mu_i=\frac{\exp(\bX_i^T\bbeta_0+f_0(t_i))}{1+\exp(\bX_i^T\bbeta_0+f_0(t_i))}.$$
The logistic link makes sense if the canonical parameter $\eta(\cdot)$ belongs to the interval $[-4,4]$, as one can see for example page 28 of \cite{FahrmeirTutz}. 
Consequently, to get a SNR of order $3$ one may choose a parameter $m$ in the binomial distribution equal approximatively to $50$. This choice is not adapted for real data applications, especially when the binomial regression corresponds to a classification problem.

Due to this remark, we choose here to make a compromise and to apply the algorithm with a parameter $m$ equal to $20$, which corresponds to $21$ classes. This choice is not coherent with a classification problem but allows to consider much reasonable SNRs. 
The observations of a simulated sample are represented in Figure~\ref{fig:obs_gplm_bin}. The results for this example are summarized in Table~\ref{tab:table_bino}. Figure~\ref{fig:boxplot_bino} gives the boxplots of the squared error obtained by the five methods considered.

\begin{figure}[!h]
\vspace{-12pt}
\centering{   \includegraphics[width=150mm,height=45mm]{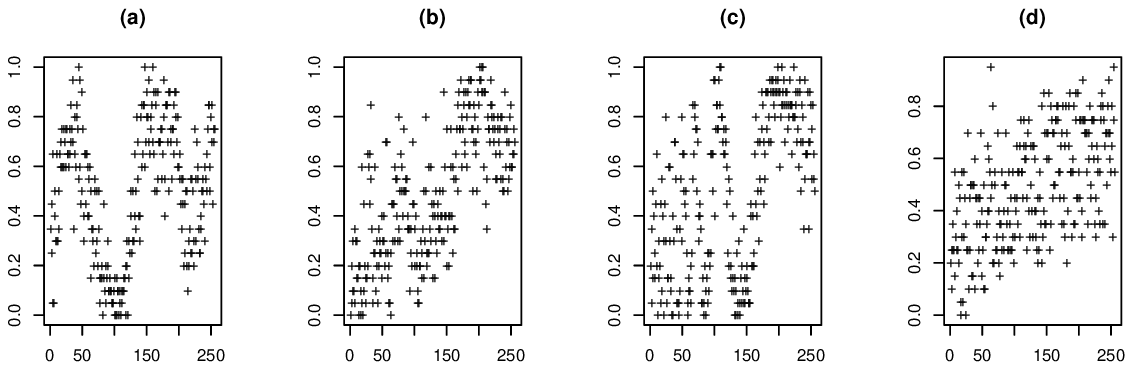}}
\caption{An example of the observations obtained on a simulation in Example 2, with the {\it Heavisine} function in Figure (a), the {\it Blocks} function in Figure (b), the {\it Doppler} function in Figure (c) and the {\it Bumps} function in Figure (d).}
\label{fig:obs_gplm_bin}       
\end{figure}

\begin{table}[!h]

\caption{Measures of quality the estimates over the 500 simulations in Example 2 with $n=2^{8}$ and $\phi=1/20$.}

\label{tab:table_bino}

\begin{center}
\begin{tabular}{cc}
\begin{minipage}{0.5\textwidth}{
\begin{tabular}{@{} l@{\hspace{0.5cm}} l @{\hspace{0.5cm}} l @{\hspace{0.5cm}}l @{}}

\hline\noalign{\smallskip}

&\multicolumn{3}{c}{{\it Heavisine}} \\ [5pt]
   & {$MSE_{\bbeta}$} & {$AMSE_f$} & Global AMSE\\
\noalign{\smallskip}\hline\noalign{\smallskip}
  {\it Wavelets} & 0.1538  & 0.1172&0.0691\\
  {\it Kernel Speckman} &0.0305 & 0.0640 & 0.0621\\
  {\it Kernel Backfitting} & 0.1134 & 0.0932 & 0.0830\\
 {\it Spline Speckman} & {\bf 0.0244} & {\bf 0.0349} &{\bf 0.0339}\\
  {\it Spline Backfitting} & 0.0294 & 0.0350 & 0.0355\\
 
 \noalign{\smallskip}\hline\\
 
&\multicolumn{3}{c}{{\it Blocks}} \\ [5pt]
   & {$MSE_{\bbeta}$} & {$AMSE_f$} & Global AMSE\\
\noalign{\smallskip}\hline\noalign{\smallskip}
  {\it Wavelets} & 0.0230 & 0.1945 & 0.1811\\
  {\it Kernel Speckman} & 0.0341 & 0.2511 &0.2511\\
  {\it Kernel Backfitting} & {\bf 0.0184} & 0.2145 &0.2126\\
 {\it Spline Speckman} & 0.0706 & 0.2096 & 0.1981\\
  {\it Spline Backfitting} & 0.0304 & {\bf 0.1314} & {\bf 0.1274}\\

 \noalign{\smallskip}\hline
\end{tabular}
}\end{minipage}
&
\begin{minipage}{0.5\textwidth}{
\begin{tabular}{@{} l@{\hspace{0.5cm}} l @{\hspace{0.5cm}} l @{\hspace{0.5cm}}l @{}}

\hline\noalign{\smallskip}

&\multicolumn{3}{c}{{\it Doppler}} \\ [5pt]
   & {$MSE_{\bbeta}$} & {$AMSE_f$} & Global AMSE\\
\noalign{\smallskip}\hline\noalign{\smallskip}
  {\it Wavelets} & {\bf 0.0159}  & 0.1801 & 0.1846\\
  {\it Kernel Speckman} & 0.0427 & 0.03662 &0.3737\\
  {\it Kernel Backfitting} & 0.0421 & 0.3486 &0.3488\\
 {\it Spline Speckman} & 0.0444 & 0.2060 &0.1914\\
  {\it Spline Backfitting} & 0.0233 & {\bf 0.1780} & 0.1762\\

 \noalign{\smallskip}\hline\\

&\multicolumn{3}{c}{{\it Bumps}} \\ [5pt]
   & {$MSE_{\bbeta}$} & {$AMSE_f$} & Global AMSE\\
\noalign{\smallskip}\hline\noalign{\smallskip}
  {\it Wavelets} & {\bf 0.0115} & 0.0510 & 0.0510\\
  {\it Kernel Speckman} & 0.0200 & 0.0492 &0.0474\\
  {\it Kernel Backfitting} &0.0170 & {\bf 0.0490} & {\bf 0.0467}\\
 {\it Spline Speckman} & 0.0153 & 0.0520 &0.0497 \\
  {\it Spline Backfitting} & 0.0165 & 0.0511 &0.0497\\

 \noalign{\smallskip}\hline

\end{tabular}}\end{minipage}
\end{tabular}
\end{center}

\end{table}


\begin{figure}[h]
\centering{
{\it Squared Error for $\bbeta$}\\
  \includegraphics[width=150mm,height=35mm]{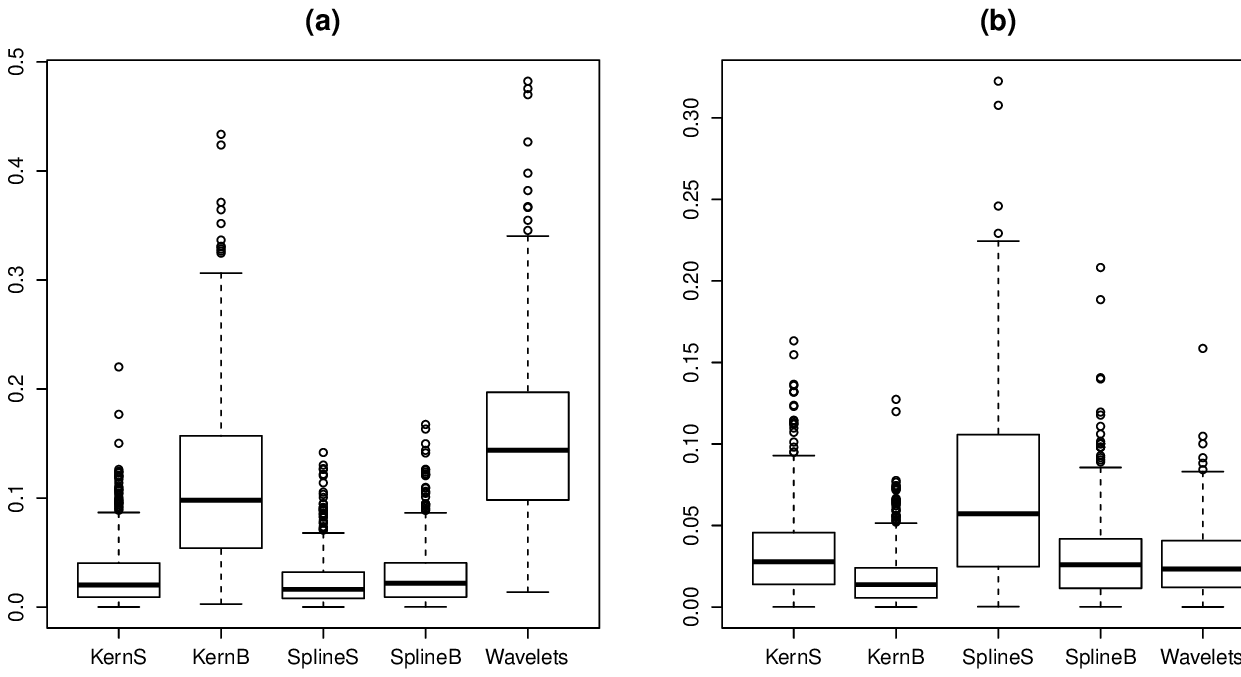}\\
{\it Mean Squared Error for the functional part }\\ 
 \includegraphics[width=150mm,height=35mm]{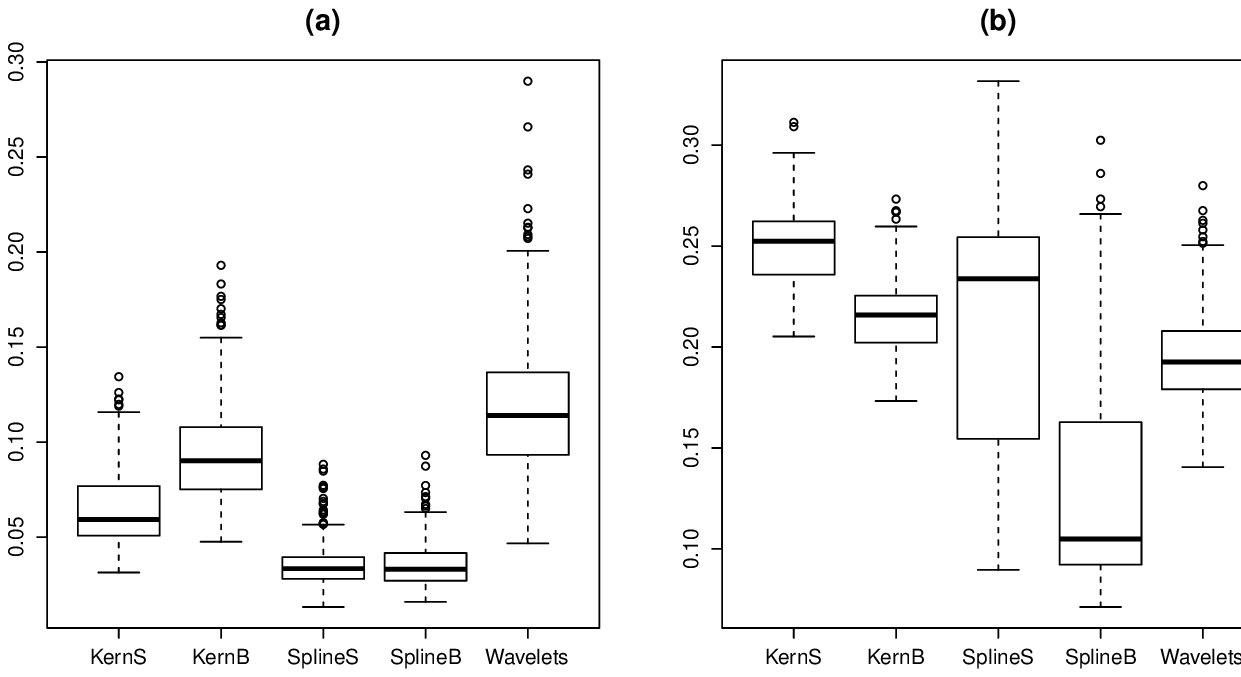}\\
{\it Global Mean Squared Error}\\
 \includegraphics[width=150mm,height=35mm]{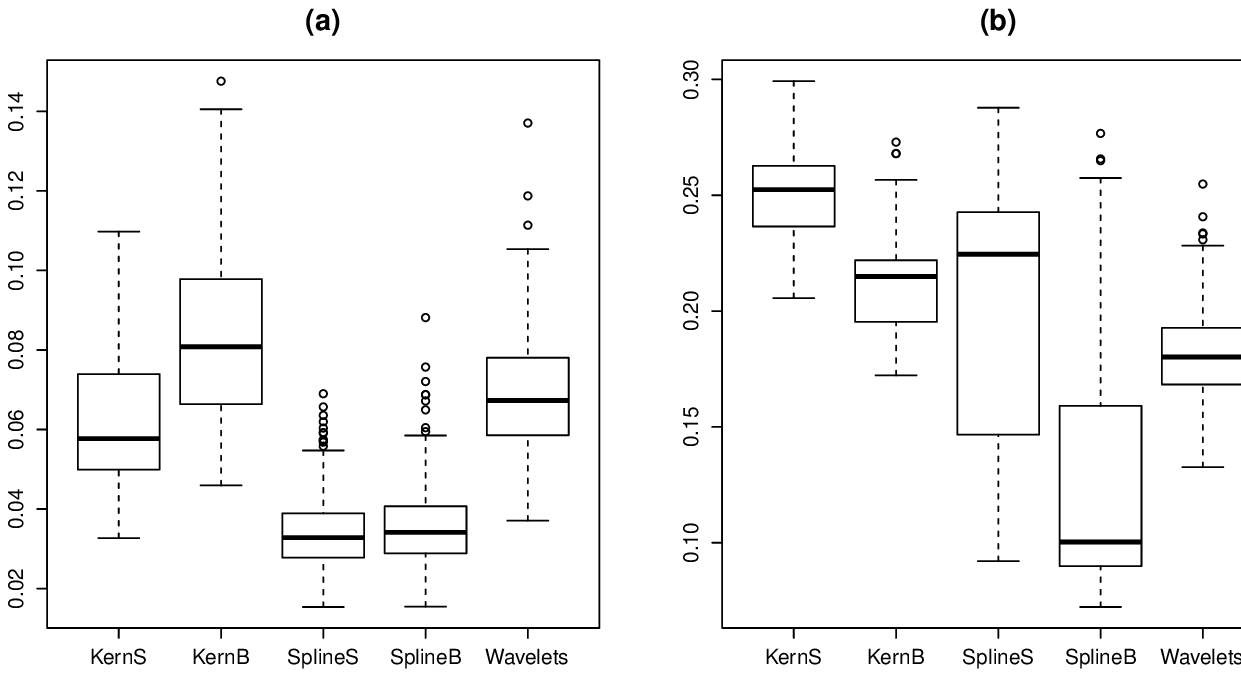}\\
}
\caption{Boxplots of the MSE in Example 2, with the {\it Heavisine} function in Figures (a), the {\it Blocks} function in Figures (b), the {\it Doppler} function in Figures (c) and the {\it Bumps} function in Figures (d). {\it KernS} and {\it KernB} stand for the kernel procedures respectively with Speckman algorithm and Backfitting algorithm, and {\it SplineS} and {\it SplineB} stands for the splines procedures respectively with Speckman algorithm and Backfitting algorithm.}
\label{fig:boxplot_bino}       
\end{figure}

\begin{figure}[!h]
\vspace{-12pt}
\centering{ \includegraphics[width=150mm,height=45mm]{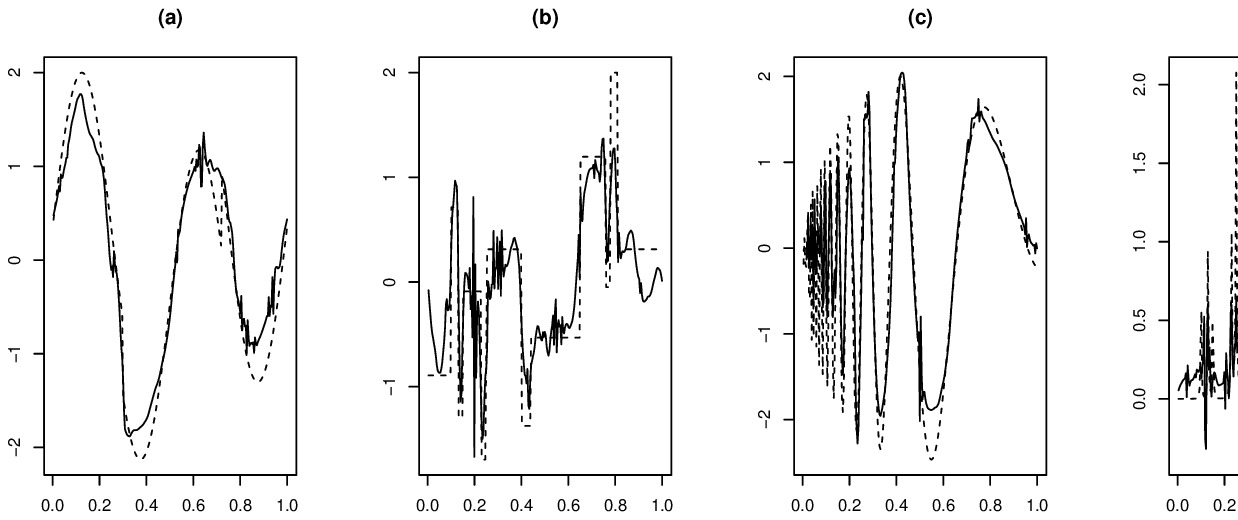}}
\caption{An example of estimation of the nonparametric part in Example 2, with the {\it Heavisine} function in Figure (a), the {\it Blocks} function in Figure (b), the {\it Doppler} function in Figure (c) and the {\it Bumps} function in Figure (d).  Dots lines corresponds to the true functions and plain lines to their estimate.}
\label{fig:gplm_bino}       
\vspace{-12pt}

\end{figure}

Contrary to Example 1 where the kernel procedures over-performed the splines in the functional part estimation, it appears in Example 2 that the splines-based procedures give better quality indexes. The case of the {\it Doppler} signal is the most relevant of this fact. On the same way, the differences which appeared in Example 1 between backfitting algorithms and Speckman algorithms are no more observed. The backfitting algorithm even give better results, except with the {\it Heavisine} signal.

Concerning the wavelets scheme, one can see that with the {\it Heavisine} functional part, it still has the worst quality for the linear and the functional part. Yet the difference with others methods is less important by comparison with the Gaussian distribution. With the others nonparametric test functions, the squared errors of the wavelets-based estimators are nearly the same as splines-based estimators and they lead to better quality indexes than kernel procedures.

Figure~\ref{fig:gplm_bino} presents an example of estimation of the functional parts for one simulation with each test functions. Contrary to the Gaussian distribution, the threshold level chosen here leads to under-smoothing of the estimated signals. Nevertheless the visual quality seems satisfying.




\subsubsection{Example 3: Poisson distribution}



In Example 3, we consider a Poisson distribution: $y_i\sim\mathcal P(\mu_i)$ with $\mu_i=\exp(\eta_i)$. Observations of a simulated data sample are represented in Figure~\ref{fig:obs_gplm_poisson}. 


\begin{figure}[!h]
\centering{  \includegraphics[width=150mm,height=45mm]{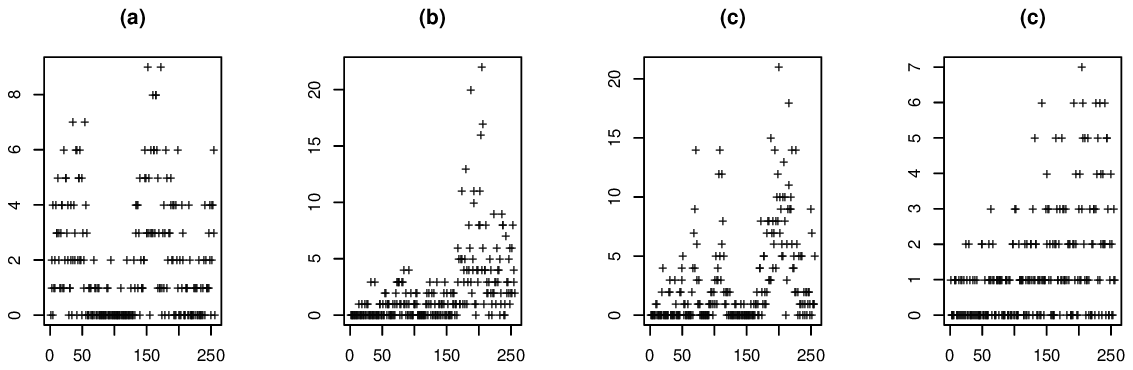}}
\caption{An example of the observations obtained on a simulation in Example 3, with the {\it Heavisine} function in Figure (a), the {\it Blocks} function in Figure (b), the {\it Doppler} function in Figure (c) and the {\it Bumps} function in Figure (d).}
\label{fig:obs_gplm_poisson}       
\end{figure}

Wavelets procedure with a Poisson distribution is highly instable. This is due to the fact that the inverse link function in this model is exponential. If an observation has a high value, the inverse link transform will explode. The wavelets procedure will detect this as a high frequency element in the functional part. Due to the numerical limitations of the software this phenomenon can return infinite values, even if they should be thresholded at the next step. As a consequence, we modify the algorithm by imposing the values returned by the inverse link and the variance functions to be reasonable. With kernel or splines based procedure this problem does not occur since they sufficiently smooth the data. 

The results for this example are summarized in Table~\ref{tab:table_poisson} and the boxplots of the squared errors indexes are given respectively for the linear regressor, the functional component and the global estimation for $\eta$ in Figure~\ref{fig:boxplot_poisson}. Even without considering the boundaries of the observation plan to avoid boundaries effects, we obtain similar values. In particular, this is the case for the kernel estimators with Speckman algorithm. 

\begin{table}[h]

\caption{Measures of quality the estimates over the 500 simulations in Example 3 with $n=2^{8}$. Lowest values are in bold face type.}

\label{tab:table_poisson}

\begin{center}
\vspace{-6pt}

\begin{tabular}{cc}
\begin{minipage}{0.5\textwidth}{
\begin{tabular}{@{} l@{\hspace{0.5cm}} l @{\hspace{0.5cm}} l @{\hspace{0.5cm}}l @{}}

\hline\noalign{\smallskip}

&\multicolumn{3}{c}{{\it Heavisine}} \\ [5pt]
   & {$MSE_{\bbeta}$} & {$AMSE_f$} & Global AMSE\\
\noalign{\smallskip}\hline\noalign{\smallskip}
  {\it Wavelets} & 1.0231  & 1.0434 & 0.6860 \\
  {\it Kernel Speckman} & 0.0592 & 1.5667 &1.5611\\
  {\it Kernel Backfitting} & 0.0352  & 0.9048 &0.8951\\
 {\it Spline Speckman} & 0.05388  & 0.3762 & 0.3793\\
  {\it Spline Backfitting} & {\bf 0.0226}  & {\bf 0.2881} & {\bf 0.2849}\\ 
 
 \noalign{\smallskip}\hline \\ 

&\multicolumn{3}{c}{{\it Blocks}} \\ [5pt]
   & {$MSE_{\bbeta}$} & {$AMSE_f$} & Global AMSE\\
\noalign{\smallskip}\hline\noalign{\smallskip}
  {\it Wavelets} & 0.4345 & 0.7637 & 0.6393\\
  {\it Kernel Speckman} & 0.2649 & 7.7513 &8.0241\\
  {\it Kernel Backfitting} & 0.0561 & 0.3581 &0.3460\\
 {\it Spline Speckman} & {\bf 0.0508 }& 0.8458 & 0.8348\\
  {\it Spline Backfitting} & 0.0580 & {\bf 0.2774} &{\bf  0.2693}\\

 \noalign{\smallskip}\hline 
 
 \end{tabular}}\end{minipage} 
 
 &
\begin{minipage}{0.5\textwidth}{
\begin{tabular}{@{} l@{\hspace{0.5cm}} l @{\hspace{0.5cm}} l @{\hspace{0.5cm}}l @{}}

\hline\noalign{\smallskip}

&\multicolumn{3}{c}{{\it Doppler}} \\ [5pt]
   & {$MSE_{\bbeta}$} & {$AMSE_f$} & Global AMSE\\
\noalign{\smallskip}\hline\noalign{\smallskip}
  {\it Wavelets} & 0.2482  & 1.6852 & 1.6124\\
  {\it Kernel Speckman} & 0.1376 & 9.4225 &9.4966\\
  {\it Kernel Backfitting} & 0.2337 & 0.1139 & 2.0540\\
 {\it Spline Speckman} & {\bf 0.0336} &{\bf  0.6057 }&{\bf 0.5944}\\
  {\it Spline Backfitting} & 0.0444 & 0.7566 & 0.7525\\

 \noalign{\smallskip}\hline \\

&\multicolumn{3}{c}{{\it Bumps}} \\ [5pt]
   & {$MSE_{\bbeta}$} & {$AMSE_f$} & Global AMSE\\
\noalign{\smallskip}\hline\noalign{\smallskip}
  {\it Wavelets} & 0.1595 & 0.3970 & 0.3812\\
  {\it Kernel Speckman} & 0.6483 & 2.3864  &2.4557\\
  {\it Kernel Backfitting} & 0.0620 & {\bf 0.1139} &{\bf 0.0998}\\
 {\it Spline Speckman} & 0.0680 & 0.1702 &0.1598\\
  {\it Spline Backfitting} & {\bf 0.0553} & 0.1461 &0.1374\\
 
\noalign{\smallskip}\hline 
 
\end{tabular}}\end{minipage}

\end{tabular}

\end{center}
\vspace{-3cm}
\end{table}

\begin{figure}[h]
\centering{
{\it Squared Error for $\bbeta$}\\
  \includegraphics[width=150mm,height=35mm]{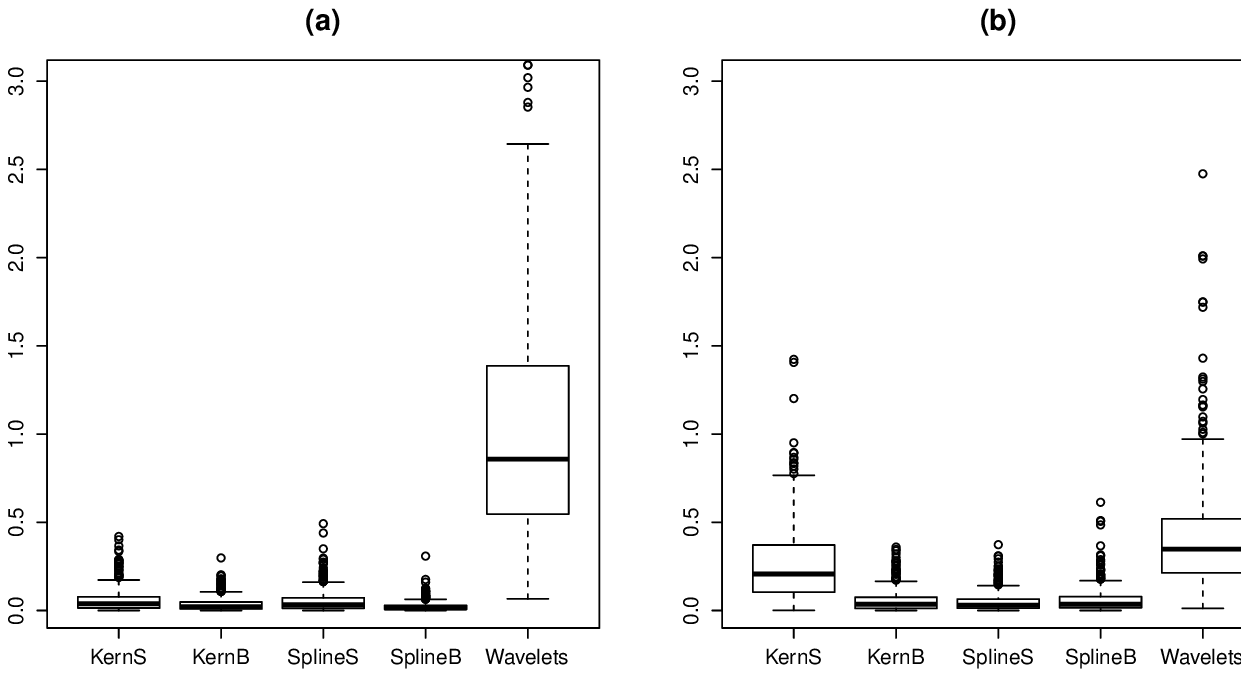}\\
{\it Mean Squared Error for the functional part }\\ 
 \includegraphics[width=150mm,height=35mm]{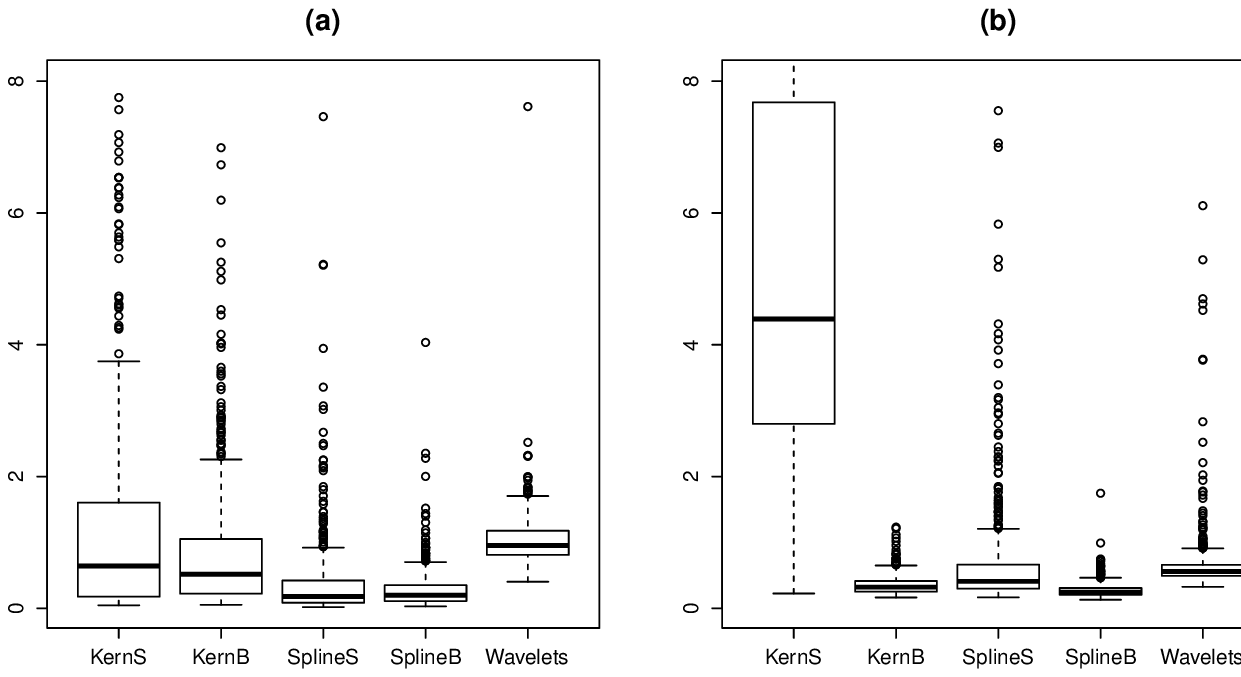}\\
{\it Global Mean Squared Error}\\
 \includegraphics[width=150mm,height=35mm]{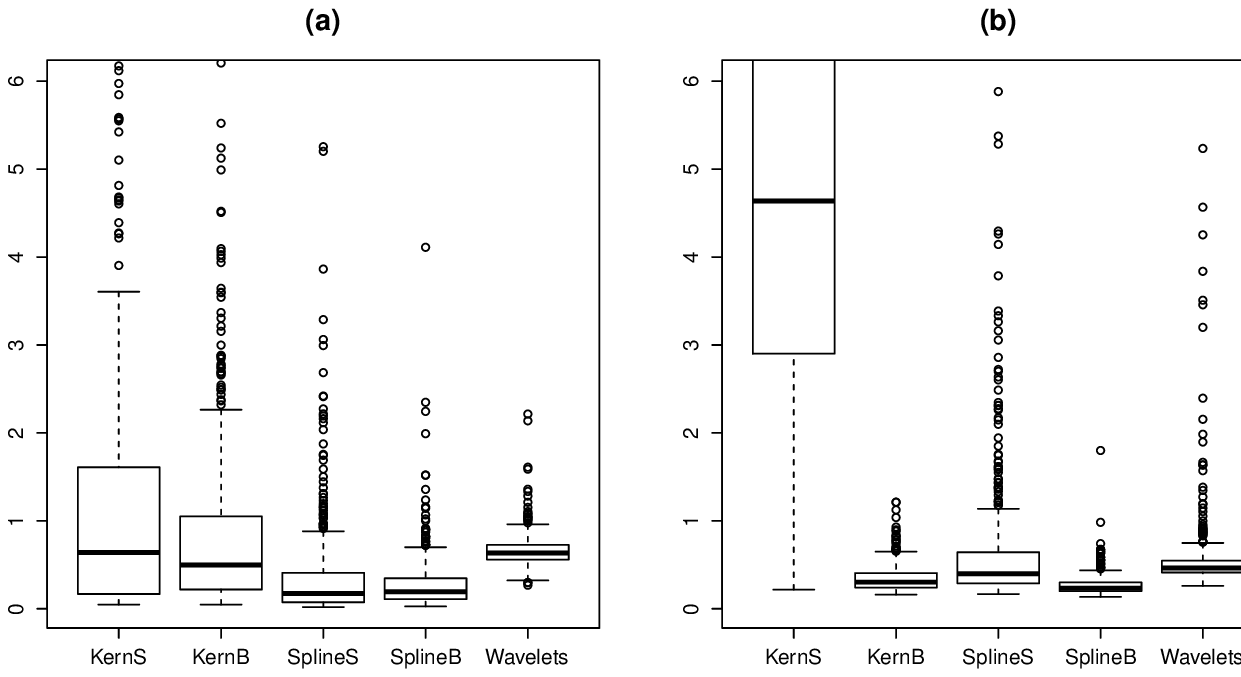}\\
}
\caption{Boxplots of the MSE in Example 3, with the {\it Heavisine} function in Figures (a), the {\it Blocks} function in Figures (b), the {\it Doppler} function in Figures (c) and the {\it Bumps} function in Figures (d). The plan [0.2, 0.8] was used in the calculation of the MSEs to avoid boundaries effects. {\it KernS} and {\it KernB} stand for the kernel procedures respectively with Speckman algorithm and Backfitting algorithm, and {\it SplineS} and {\it SplineB} stands for the splines procedures respectively with Speckman algorithm and Backfitting algorithm.}
\label{fig:boxplot_poisson} 
\end{figure}

\begin{figure}[!h]
\centering{ \includegraphics[width=150mm,height=45mm]{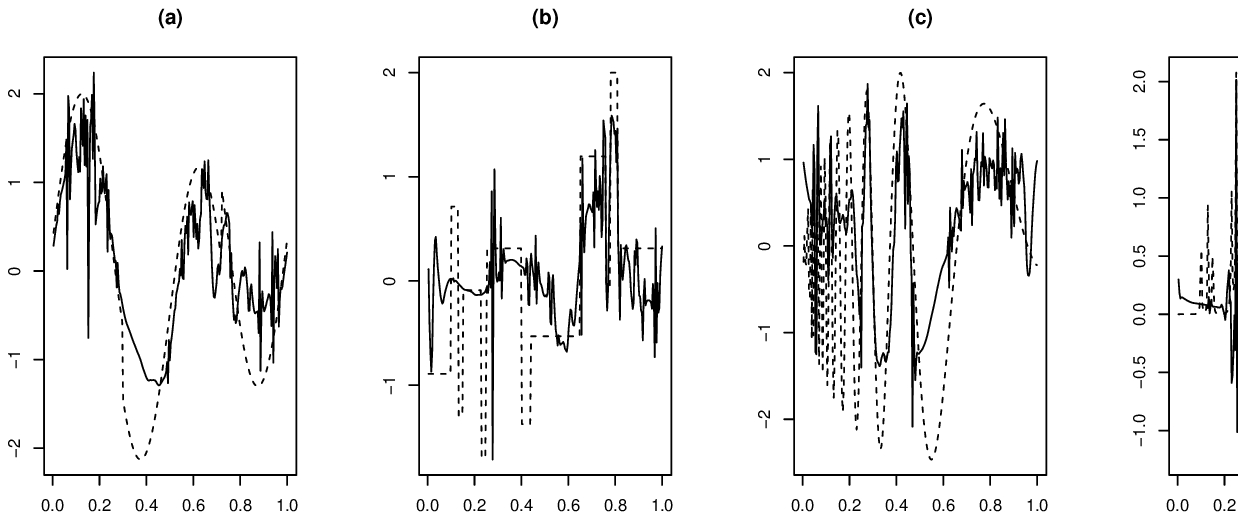}}
\caption{An example of estimation of the nonparametric part in Example 3, with the {\it Heavisine} function in Figure (a), the {\it Blocks} function in Figure (b), the {\it Doppler} function in Figure (c) and the {\it Bumps} function in Figure (d).}
\label{fig:gplm_poisson}       
\end{figure}

\newpage

It is possible that the very bad quality indexes for the kernel estimators with Speckman algorithm are due to non robustness and could be avoid by forcing the inverse link function to be reasonable, as it was done for wavelets estimators. Yet this was not implemented. The fact we artificially bound the inverse link function can moreover introduce bias.

Wavelets based estimation in the Poisson model are less performant than splines ones, as it can be seen in Table~\ref{tab:table_poisson} and in Figure~\ref{fig:boxplot_poisson}. In particular the estimation of the linear part when the functional part is the {\it Heavisine} signal has a high Mean Squared Error. 
For the {\it Doppler} and the {\it Bumps} test functions, the wavelets estimates have a quality more comparable with splines based estimators, even if lower. And comparing to kernel procedures, especially the one with Speckman algorithm, they give more stable results. Nevertheless splines based estimators and overall with the backfitting algorithm lead to the best quality indexes.

Figure~\ref{fig:gplm_poisson} gives an example of the estimation of the functional part obtained. The visual approximation is not very good and if the global forms of the functions are retrieved, some peaks could lead to bad interpretations in a real data application

Contrary to the Gaussian and the Binomial cases where the wavelets estimators quality were satisfying, it seems much more lower with the Poisson distribution. To study this phenomenon, we simulate Generalized Functional Models with Poisson distributions, meaning that we omit the linear component in the model. In Table~\ref{tab:gfm_poisson} one can see that the quality indexes in the functional estimation are much more better than in the partially linear framework. This is confirmed with Figure~\ref{fig:boxplot_gfm_poisson} where the comparison of MSEs boxplots between GFM and GPLM shows an important difference in favour of the functional model. Moreover for readability reasons, we take truncated vertical axes. Actually, in a GPLM, wavelets estimators have higher MSE than the axes limit considered, on the contrary to the MSE in functional frameworks. Finally, visually the estimation in a functional framework is also better, comparing Figure~\ref{fig:gplm_poisson} and~\ref{fig:gfm_poisson}.

 \begin{table}[h]

\caption{Numerical indexes for the regression of the threshold level that minimizes the functional $AMSE$ with respect to $\sqrt{\log(n)}$ with a Poisson distribution. Calculations were done on $100$ samples for each size of $n$.}

\label{tab:gfm_poisson}
\begin{center}

\hspace{0.5cm}

\begin{tabular}{@{} l@{\hspace{0.5cm}} llll @{}}
\hline\noalign{\smallskip}
  Function & {\it Heavisine} & {\it Blocks} & {\it Doppler} & {\it Bumps}\\
  AMSE in a GFM & 0.4368 & 0.3730 & 0.5306 & 0.1549 \\
 AMSE in a GPLM & 1.0434 & 0.7637 & 1.6852 &  0.3970\\
\noalign{\smallskip}\hline
\end{tabular}
\end{center}
\end{table}

\begin{figure}[h]
\centering{
 \includegraphics[width=150mm,height=40mm]{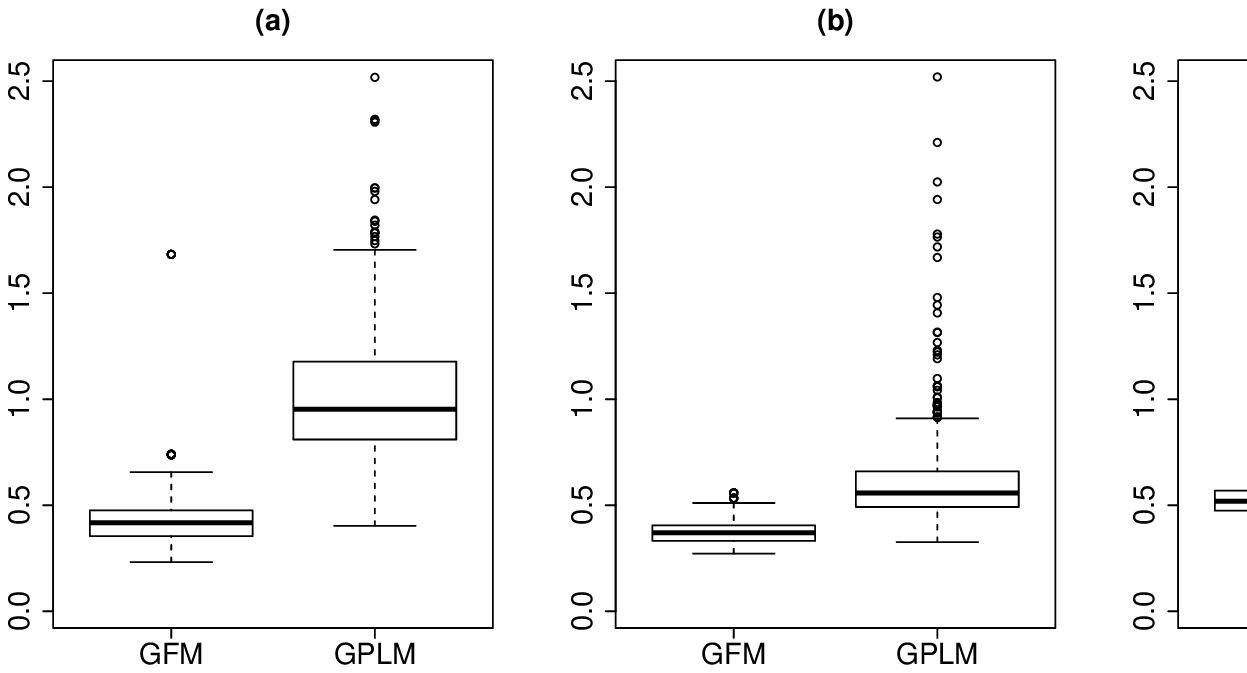}\\

}
\caption{Boxplots of the functional MSE for a Poisson distribution, in a GFM and in a GPLM, with the {\it Heavisine} function in Figures (a), the {\it Blocks} function in Figures (b), the {\it Doppler} function in Figures (c) and the {\it Bumps} function in Figures (d).}
\label{fig:boxplot_gfm_poisson}       
\end{figure}

\begin{figure}[!h]
\centering{    \includegraphics[width=150mm,height=45mm]{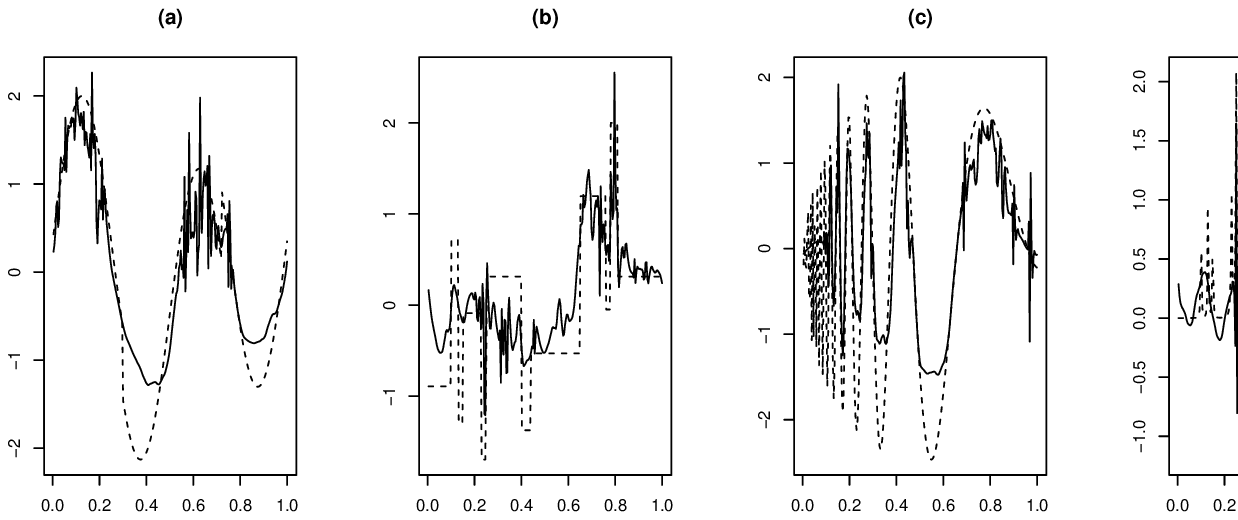}}
\caption{An example of estimation of the nonparametric part in a Poisson functional regresion, with the {\it Heavisine} function in Figure (a), the {\it Blocks} function in Figure (b), the {\it Doppler} function in Figure (c) and the {\it Bumps} function in Figure (d).}
\label{fig:gfm_poisson}       
\end{figure}

A possible interpretation of this phenomenon is the lack of robustness of the wavelets. As explained previously a bound in the inverse link function is necessary for numerical implementation in Poisson GPLM. This is due to the fact that high values data is considered as high frequencies. This artificial bound can imply a bias, but is relevant overall of a non robustness in the procedure. Moreover, the high values obtained in the functional estimation will be considered as outliers in the linear part estimation step (see \cite{gannaz}) and imply a bias estimation of the linear regressor. Consequently, a more robust procedure as suggested by \cite{BoenteHeZou} could improve the proposed estimation scheme.








\newpage

\section*{Conclusion}

This paper proposes a penalized loglikelihood estimation in generalized partially linear models. With appropriate penalties, we establish the asymptotic near-minimaxity of the estimators for both the linear and the nonparametric parts of the model. The result holds for a large class of functions, possibly non smooth and the conditions of correlation between the covariate design of the linear part and the functional part are similar to literature. Moreover with an $\ell^1$ penalty on the wavelet coefficients of the function (leading to soft thresholding), the procedure appears to be adaptive relatively to the smoothness of the function. In this particular case we develop an implementation algorithm. Thanks to a simulations study we see that the performance of the estimator is satisfactory with Gaussian and Binomial distributions, even if, by comparison, kernel and splines based estimators with a Generalized Cross Validation procedure sometimes give better quality indexes. The proposed wavelets procedure also seems to give satisfactory results in a Poisson functional regression, but leads to worse quality estimates than splines procedures in a Poisson partially linear regression.

Our ongoing research may deal with a more robust estimation procedure in Poisson partially linear models, following {\it e.g.} \cite{BoenteHeZou}. The estimation of the dispersion parameter in generalized partially linear models and the study of the efficiency of the estimation of the regression vector~$\bbeta$ are also interesting perspectives. For instance, \cite{diff_finies} explored this issue for partially linear models with gaussian distributions. An extension to generalized frameworks would be interesting. Developments in non-equidistant designs for the nonparametric part should also be explored.

{\bf Aknowledgements.}
The author would like to thank the referees for their constructive criticism which helped improving the paper. Pr. Antoniadis is gratefully acknowledged for his constructive comments and fruitful discussions. The author would also like to thank Pr. Gadat for giving some helpful references.

{
\bibliographystyle{spbasic}      
\bibliography{biblio}

}


\section{Appendix}

\appendix

The proof is structured as follows: in Section A, we justify the fact that under a structure $(A_{corr})$ of the covariates, with a well-chosen penalty, the scale part of the linear part does not intervene in the quality of the estimators. Section B studies the behaviour of the nonparametric part $f$, while Section $C$ presents the asymptotic properties of the estimator for the linear regressor $\beta$.

\section{Decorrelation of the scale components}

We suppose here that the assumption $(A_{corr})$ holds and that in the penalty only intervene the wavelets coefficients. To be more precise, given a vector of real values $\be=(e_1,\ldots, e_n)^T$, the discrete wavelet transform of $\be$ is given by a $n\times1$ vector comprising both discrete scaling coefficients, noted $(\theta^S_{j_S,k})_{k\in\Z}$, and discrete wavelet coefficients, $(\theta^W_{j,k})_{j\geq j_S, \,k\in\Z}$. The wavelet inverse transform of the scaling coefficients $(\theta^S_{j_S,k})_{k\in\Z}$ (fixing the wavelet coefficients equal to zeros) is noted $e^S$ and the wavelet inverse transform of the wavelets coefficients $(\theta^W_{j_S,k})_{k\in\Z}$ (fixing the scale coefficients equal to zeros) is noted $e^W$. The penalty is in this section assumed to be such that $Pen(f)$ is a function of $f^W$ only. To better apprehend the mechanisms, let us decompose each components $\bX_i$ and $f$ in the criterion into their scale parts $\bX_i^S$ and $f^S$ and their wavelets parts $\bX_i^W$ and $f^W$

Recall that the estimation procedure consists of maximizing the criterion given in equation~(\ref{critere}), first with respect to the functional part $f$ and afterwards with respect to the linear component. For any function $f$ and any $p$-dimensional vector $\bbeta$, the criterion $K_{n,\lambda}(f,\bbeta)$ given in (\ref{critere}) can be written 
\begin{equation}\label{critere2}
 K_{n,\lambda}(f,\bbeta)=\sum_{i=1}^n \mathcal L\left(y_i, \bX_i^{WT}\bbeta+a(f,X^S\bbeta)(t_i)\right) \, - \, \lambda\,Pen(a(f,X))=:L_{n,\lambda}(a(f,X^S\bbeta),X^W\bbeta),
\end{equation}
with $a(f,X^S\bbeta)(t_i)=\bX_i^{ST}\bbeta+f(t_i).$ Note that the criterion $L_{n,\lambda}(h,X^W\bbeta)$ does not depend of $X^S$ for a given function $h$. For a fixed vector $\bbeta$, $\widetilde f_{n,\bbeta}=\argmax_{f} K_{n,\lambda}(f,\bbeta)=h_{\bbeta}-X^{ST}\bbeta$ with $h_{\bbeta}=\argmax_h L_{n,\lambda}(h,\bbeta)$. Replacing $f$ by its estimation, $K_{n,\lambda}(\widetilde f_{n,\bbeta},\bbeta)=\sum_{i=1}^n \mathcal L\left(y_i, \bX_i^{WT}\bbeta+h_{\bbeta}(t_i)\right) \, - \, Pen(h_{\bbeta})$. Consequently, the estimate $\hat\bbeta_n=\argmax_{\bbeta} K_{n,\lambda}(\widetilde f_{n,\bbeta},\bbeta)$ is in fact defined independently from $X^S$ in such a framework.

It is clear that the properties established for the argument maximizing the criterion $K_{n,\lambda}$ are available when maximizing the criterion $L_{n,\lambda}$. When the structure of the covariates and of the penalty are well-adapted we will then consider without loss of generality that the covariates $(\bX_i^W)_{i=1,\dots,n}$ contain only the wavelets components $(\bX_i^W)_{i=1,\dots,n}$.

\section[Estimation of the nonparametric part]{Estimation of the nonparametric part by $\widetilde f_{n,\bbeta}$}

In this Section, we have to distinguish between the non-adaptive type penalties of Section~\ref{sec:nonadapt} and the $\ell^1$ penalty on wavelet coefficients. The whole scheme is the same, so the second case will be less detailed.

\subsection{Nonadaptive case}

Here, the penalty is supposed directly linked with an entropy control through assumptions~(A5) and~(A6).

\subsubsection{A Taylor expansion}

The difference between the criterions $K_{n,\lambda}(f,\bbeta)$ and $K_{n,\lambda}(f_0,\bbeta)$ is equal to: $$\frac{1}{n}\left(K_{n,\lambda}(f,\bbeta)-K_{n,\lambda}(f_0,\bbeta)\right)= \frac{1}{n}\sum_{i=1}^n \mathcal L(y_i,\bX_i^T\bbeta+f(t_i))-\sum_{i=1}^n \mathcal L(y_i,\bX_i^T\bbeta_0+
f_0(t_i))-\frac{\lambda}{n}\left( Pen(f)-Pen(f_0)\right).$$ Note $\eta_{0,i}=\bX_i^T\bbeta_0+f_0(t_i)$.
A Taylor expansion of degree 2 at the points $(\mathcal L(y_i,\eta_{0,i}))_i$ gives:
\begin{eqnarray*}
\frac{1}{n}\left(K_{n,\lambda}(f,\bbeta)-K_{n,\lambda}(f_0,\bbeta)\right) \, = & &\underbrace{\frac{1}{n}\sum_{i=1}^n \dot{\mathcal L}(y_i,\eta_{0,i})\left(f(t_i)-f_0(t_i)\right)}_{A_1} - \underbrace{1/2\frac{1}{n}\sum_{i=1}^n \ddot b(\eta_{0,i})\left(f(t_i)-f_0(t_i))\right)^2}_{A_2} \\
& - &\underbrace{\frac{1}{n}\sum_{i=1}^n \ddot b(\eta_{0,i})\left(\bX_i^T(\bbeta-\bbeta_0)\right)\left(f(t_i)-f_0(t_i)\right)}_{Co} - \frac{\lambda}{n} \left(Pen(f)-Pen(f_0)\right) + T_1(f,\bbeta)
\end{eqnarray*}
with $T_1(f,\bbeta)\,\leq \,\frac 1 6 \sup_i\sup_\eta \dddot b(\eta_{i})\, \frac 1 n \sum_i\left( \left|\bX_i^T(\bbeta-\bbeta_0)+f(t_i)-f_0(t_i)\right|^3+\left|f(t_i)-f_0(t_i)\right|^3\right).$

Let us study the term $T_1(f,\bbeta)$. Under assumption (A3), $T_1(f,\bbeta)\leq cst \frac 1 n \sum \left(\left|\bX_i^T(\bbeta-\bbeta_0) + f(t_i)-f_0(t_i)\right|^3+\left|\bX_i^T(\bbeta-\bbeta_0)\right|^3\right).$
Using the convexity of $x\mapsto |x|^3$, it comes:
$$T_1(f,\bbeta)\leq cst \left(\underbrace{\frac 1 n \sum_{i=1}^n \left|\bX_i^T(\bbeta-\bbeta_0)\right|^3}_{T_{11}}+\underbrace{\frac 1 n \sum_{i=1}^n \left|f(t_i)-f_0(t_i)\right|^3}_{T_{12}}\right).$$
Then, one can easily see that when assumptions (A1) and (A2) hold, $T_{11}=o(ph^{1/2}/n)$.
Consequently, if $ph^{1/2}v_n^2/n$ is bounded, then for every function $f$ satisfying $\frac{v_n^2}{n}\sum_i\left|f(t_i)-f_0(t_i)\right|^3\leq c$, we have $v_n^2 T_1(f,\bbeta)\leq C.$
When $v_n=n^{1/(2+\nu)}$ with $\nu>0$, the sequence $ph^{1/2}v_n^2/n$ is bounded if $ph^{1/2}n^{-\nu/(2+\nu)}$ is bounded.

Let $$S_n(f)=\frac{v_n^2}{n}\left(K_{n,\lambda}(f,\bbeta)-K_{n,\lambda}(f_0,\bbeta)-A_1+A_2-Co\right)+\frac{\lambda v_n^2}{n} \left(Pen(f)-Pen(f_0)\right).$$ For every sequence $u_n\to 0$ and for all function $f$ satisfying $\frac{v_n^2}{n}\sum_i\left|f(t_i)-f_0(t_i)\right|^3\leq c$, we prove $u_n S_n(f)\to 0$. Fix $u_n$ a decreasing sequence going to 0. Convexity argumentation (Theorem 10.8 in \cite{Rockafellar} with a change of variable on $f$ to obtain a set independent from $n$) gives: $\sup_{\frac{v_n^2}{n}\sum_i\left|f(t_i)-f_0(t_i)\right|^3\leq c} u_n S_n(f)\to 0$. Note now that $\frac{v_n^2}{n}\sum_i\left|f(t_i)-f_0(t_i)\right|^3\leq v_n^2\|f-f_0\|_n^2 \|f-f_0\|_\infty,$ and that we consider the maximisation problem on the set $\{f,\,\|f\|_\infty\leq C_\infty\}$.
\begin{equation}\label{conv_sup}\sup_{\{f,\,N(f)\leq c\,\|f\|_\infty\leq C_\infty} u_n S_n(f)\to 0, \text{~ where~} N(f)=v_n^2\|f-f_0\|_n^2+\ddot b_0^{-1} J(f),\end{equation}
where $\ddot b_0=\inf_{i=1,\dots,n}\ddot b(\eta_{0,i})$.

{\bf Remark:} We need $f\mapsto S_n(f/v_n)$ to be concave for the norm $\|.-f_0\|_n$ for every $n$. It is worthy noticing the penalty $Pen$ does not intervene in the definition of $S_n$.

\subsubsection{Behaviour of the different terms}

Similarly to \cite{BaiRaoWu}, we suppose $N(f)\geq c_n$ with $c_n u_n\to \infty$.
Then one can build a sequence $c'_n$ satisfying $c'_n u_n\to\infty$ et $c'_n\leq c_n$ and such that $\sup_{N(f)\leq c'_k} u_k S_k$ goes to 0 when $k$ goes to infinity.

For $N(f)= c_n'$, this means: \begin{equation}{v_n^2}\left(\frac{1}{n}\left(K_{n,\lambda}(f,\bbeta)-K_{n,\lambda}(f_0,\bbeta)\right)-A_1+A_2-Co\right)+\frac{\lambda v_n^2}{n} \left(Pen(f)-Pen(f_0)\right)=\bigcirc(1).\end{equation} 

\begin{itemize}
\item{\bf Control of the term $A_1=\frac{1}{n}\sum_{i=1}^n \dot{\mathcal L}(y_i,\eta_{0,i})\left(f(t_i)-f_0(t_i)\right)$.}\\
Following \cite{MammenVanDerGeer}, let $\mathcal A$ be a set of uniformly bounded functions on $[0,1]$ satisfying $$\lim\sup_{n\to\infty}\,\sup_{\delta>0}\,\delta^\nu\mathcal H(\delta,\mathcal A,\|.\|_n)<\infty.$$ Suppose $\frac{f}{1+J(f)}\in \mathcal A$. Then, similarly to \cite{MammenVanDerGeer}, when assumption (A4.1) holds:
$$\sup_{f,\|f\|_\infty\leq C_\infty} \sqrt{n}\frac{A_1}{\left(\frac{\|f-f_0\|_n}{1+J(f)}\vee n^{-1/(2+\nu)}\right)^{1-\nu/2}}\,=\,\bigcirc_\P(1).$$
Under assumption (A6)
, if $N(f)=c'_n$, then $J(f)\leq c'_n$, and consequently, we obtain that with $v_n=n^{1/(2+\nu)}$: $$\sup_{\{f,\|f\|_\infty\leq C_\infty,\, N(f)\leq c'_n\}} v_n^2 A_1\,=\,o_\P(c'_n).$$

\item{\bf Control of the term $A_2=1/2\frac{1}{n}\sum_{i=1}^n \ddot b(\eta_{0,i})\left(f(t_i)-f_0(t_i)\right)^2.$}\\
We immediately get: $v_n^2 A_2\geq \ddot b_0 \left(v_n\|f-f_0\|_n\right)^2$. And thus~~ $\inf_{N(f)= c'_n} v_n^2 A_2+\frac{\lambda v_n^2}{n} \left(Pen(f)-Pen(f_0)\right)\geq \ddot b_0 {c'_n}^2.$ 

\vspace{10pt}

\item{\bf Control of the term $C_0=\frac{1}{n}\sum_{i=1}^n \ddot b(\eta_{0,i})\left((\bbeta-\bbeta_0)^T\bX_i\right)\left(f(t_i)-f_0(t_i)\right).$}\\
Cauchy-Schwarz inequality applied to $Co$ gives $C_0\leq \left\|f-f_0\right\|_n \left(\frac 1 n  \sum\ddot b(\eta_{0,i}) \|\bX_i\|^2\right)^{1/2}\|\bbeta-\bbeta_0\|.$ 
When $\sqrt{n}\|\bbeta-\bbeta_0\|\leq c$ and assumption (A1) holds, we obtain that $C_0\leq C(\bbeta) p^{1/2} n^{-1/2}\left\|f-f_0\right\|_n,$ with $C(\bbeta)$ independent from $f$ and $p$. Thus: $$\sup_{N(f)\leq c'_n} v_n^2 C_0=\bigcirc_{p.s.}\left({c_n'}^{1/2}\right)=o_\P(c'n),$$ for $v_n=n^{\frac{1}{2+\nu}}$, provided $p n^{-\nu/(2+\nu)}$ is bounded.
\end{itemize}

\subsubsection{Conclusion}

Using together the convergence (\ref{conv_sup}), the bounds of $A_1$, $A_2$ and $Co$ and assumption (A5), we obtain: $$\sup_{\left\{f,\,N(f)=c'_n\right\}} u_n\frac{v_n^2}{n}\left(K_{n,\lambda}(f,\bbeta)-K_{n,\lambda}(f_0,\bbeta)\right)<0$$ with a probability going to 1 when $n$ goes to infinity. Concavity of assumption~(A6) (because it implies the decrease of the slopes) allows to extend this result to $$\sup_{\left\{f,\,N(f)\geq c'_n\right\}} u_n\frac{v_n^2}{n}\left(K_{n,\lambda}(f,\bbeta)-K_{n,\lambda}(f_0,\bbeta)\right)<0.$$ The estimator $\widetilde f_{n,\bbeta}$ is the argument realizing the maximum of $K_{n,\lambda}(\cdot,\bbeta)$ and so $\P\left(N(\widetilde f_{n,\bbeta})\geq c_n'\right)\to 0.$

\subsection[Adaptive case]{Adaptive case, with an $\ell^1$ penalty}

In this section $Pen(f)=\sum_{i=i_s}^{n}|\theta_i^W|=\sum_{j=j_S}^{J-1}\sum_{k=0}^{2^j-1} |\theta_{j,k}^W|$ with $\btheta^W$ the vector of the wavelet coefficients of $f$.

\subsubsection{Study of the penalty}

In the first time, we are willing to study how this penalty behave for the target function $f_0$. The underlying idea is to distinguish the behaviour of the penalty among the resolution degree of the wavelets coefficients. We introduce $i_W=\left(\frac{n}{\log(n)}\right)^{1/(1+2s)}.$ If the resolution level is higher we will establish a link with the Besov semi-norm, but if the resolution level is smaller, then we are going to see that the norm $\|.\|_n$ offers a sufficient control.

\begin{itemize}
\item[$\bullet$] {\bf Lowest resolution levels} 

We are interested in the study of $\frac{\lambda v_n^2}{n}\sum_{i=i_s}^{i_W} \left(|\theta_{i}^W|-|\theta_{0,i}^W|\right)$.
Using Cauchy-Schwarz inequality, $$\left|\frac{\lambda v_n^2}{n}\sum_{i=i_s}^{i_W} \left(|\theta_{i}^W|-|\theta_{0,i}^W|\right)\right|\leq i_W^{1/2} \left(\sum|\theta_i^W-\theta_{0,i}^W|^2\right)^{1/2}= (n i_W)^{1/2} \|f-f_0\|_n.$$ It is sufficient that $i_W\leq \frac{n}{\lambda^2 v_n^2}$ to deduce $\left|\frac{\lambda v_n^2}{n}\sum_{i=i_s}^{i_W} \left(|\theta_{i}^W|-|\theta_{0,i}^W|\right)\right|\leq v_n \|f-f_0\|_n.$

\item[$\bullet$] {\bf Highest resolution levels}

We are willing to bound $\frac{\lambda v_n^2}{n}\sum_{i=i_W}^n |\theta_{0,i}^W|$. Using H\"older inequality, $\sum_k|\theta_{0,j,k}^W|\leq2^{j(1-1/\pi)}\left(\sum_k|\theta_{0,j,k}^W|^\pi\right)^{1/\pi}$ because as $\psi$ and $f$ admits compact supports the number of non-zero coefficients at a resolution level $j$ is equivalent to $2^j$. Thus $$\frac{\lambda v_n^2}{n}\sum_{i=i_W}^n |\theta_{0,i}^W|\leq \frac{\lambda v_n^2}{n^{1/2}} i_W^{-(s-1/2)} \|f_0\|_{s,\pi,\infty}.$$

Taking $i_W\leq \frac{n}{\lambda^2 v_n^2}$, and $v_n=\left(\frac{n}{\lambda^2}\right)^{s/(1+2s)}$, we obtain  $\frac{\lambda v_n^2}{n}\sum_{i=i_W}^n |\theta_{0,i}^W|\leq \|f_0\|_{s,\pi,\infty}$.

\item[$\bullet$] {\bf Conclusion}

If $a$ and $b$ are two given reals, one has $0\leq |a-b|+|b|-|a| \leq 2 \min(|b|,|a-b|)$. Applying this inequality to the penalty and using previous results, for $v_n=\left(\frac{n}{\lambda^2}\right)^{s/(1+2s)}$, $$\left|\frac{\lambda v_n^2}{n}\sum_{i=i_s}^{n} \left(|\theta_{i}^W|-|\theta_{0,i}^W|\right) - \frac{\lambda v_n^2}{n}\sum_{i=i_s}^{n} \left|\theta_{i}^W-\theta_{0,i}^W\right| \right|\leq 2 v_n \|f-f_0\|_n + 2\|f_0\|_{s,\pi,\infty}.$$

\end{itemize}

\subsubsection[Bias term]{Bias term $A_1=\frac{1}{n}\sum_{i=1}^n \dot{\mathcal L}(y_i,\eta_{0,i})\left(f(t_i)-f_0(t_i)\right)$.}

Let $N(\cdot)$ and $J(\cdot)$ be defined as follows: $N(f)=v_n^2\|f-f_0\|_n^2+ \ddot b_0^{-1} \frac{\lambda v_n^2}{n}\sum_{i=i_s}^{n} \left|\theta_{i}^W-\theta_{0,i}^W\right| $ and $J(f)=n^{-1/2}\sum_{i=i_S}^{n}|\theta_i^W|$.

We will consider the functional set $\mathcal A=\{g, J(g)\leq J_0\}$. Lemme 4.3. of \cite{LoubesVanderGeer} states that the entropy of $\mathcal A$ satisfies $\mathcal H(\mathcal A,\delta,\|.\|_n)\leq A \delta^{-2}\left(\log(n)+\log(1/\delta)\right)$ where $A$ is a constant. Consequently for all $1/n<\gamma<R$, $$\int_{\gamma}^R \mathcal H(\mathcal A,u,\|.\|_n)^{1/2}du\leq A\sqrt{2\log(n)} (\log(R)-log(\gamma)).$$

Corollary 8.3 of \cite{VanderGeer00} implies that under assumption (A4.2), there exists a constant $C_0$ depending of constants $K$ and $\sigma_O^2$ of assumption (A4.2) such that for $R>8\gamma>1/n$, and $t>2\sqrt{2}C_0\max(\log(R/\gamma),R)$, we have: \begin{equation}\label{eqn:VdG}\P\left(\sup_{\{g, g\in\mathcal A, \|g\|_n\leq R\} } \frac{1}{\sqrt{n\log(n)}}\sum_{i=1}^n\dot{\mathcal L}(y_i,\eta_{0,i})g(t_i)\,\geq t \right)\,\leq \,C_0\exp\left(-\frac{t^2}{4 C_0^2 R^2}\log(n)\right).\end{equation}

Set $g=\ddot b_0^{-1}\frac{\lambda v_n^2}{\sqrt{n}} \frac{f-f_0}{N(f)}$. According to previous developments, $g$ belongs to $\mathcal A$ and $\|g\|_n\leq 1$.
Thus, applying inequality~(\ref{eqn:VdG}) with $R=1$ and $\gamma<1/8$ constant, we can deduce that for all $t\geq 2\sqrt{2}\log(1/\gamma)C_0$, we have:  $$\P\left(\sup_{\{f, N(f)=c'_n\}} \ddot b_0^-1 \frac{\lambda}{\sqrt{\log(n)}}\frac{v_n^2 A_1}{c'_n}\,\geq t \right)\,\leq\, C_0\exp\left(-\frac{t^2}{4 C_0^2 R^2}\log(n)\right).$$
Taking $\lambda=c_0\sqrt{\log(n)}$ with $c_0>2\sqrt{2}\log(8)C_0$, we obtain that $\P\left(\sup_{\{f, N(f)=c'_n\}} v_n^2 A_1\,\geq \ddot b_0 c_n' \right)$ tends to 0 when $n$ goes to infinity.

\subsubsection{End of the proof}

The scheme is very similar to what has been presented in the nonadaptive case and will not be detailed here.

\section[Estimation of the linear part]{Estimation of the linear part: study of the behaviour of $\hat\bbeta_n$}

With a Taylor expansion of order 2,

\vspace{10pt}
$K_{n,\lambda}(\widetilde f_{n,\bbeta},\bbeta)-K_{n,\lambda}(\widetilde f_{n,\bbeta},\bbeta_0)$
\begin{eqnarray*}
 = & & \sum_{i=1}^n \mathcal L(y_i,\bX_i^{T}\bbeta+\widetilde f_{n,\bbeta}(t_i))-\sum_{i=1}^n \mathcal L(y_i,\bX_i^{T}\bbeta_0+\widetilde f_{n,\bbeta}(t_i))\\
= & &\underbrace{\sum_{i=1}^n \dot{\mathcal L}(y_i,\eta_{0,i})\left(\bX_i^T(\bbeta-\bbeta_0)\right)}_{B_1} - \underbrace{1/2\sum_{i=1}^n \ddot b(\eta_{0,i})\left(\bX_i^T(\bbeta-\bbeta_0)\right)^2}_{B_2}  - \underbrace{\sum_{i=1}^n \ddot b(\eta_{0,i})\left(\bX_i^T(\bbeta-\bbeta_0)\right)\left(\widetilde f_{n,\bbeta}(t_i)-f_0(t_i)\right)}_{\widetilde Co} +   T_2(\bbeta)
\end{eqnarray*}

\subsection{Behaviour of the different terms}
We distinguish the terms where the functional part intervenes from the terms where only the linear part appears.

\subsubsection[Control of the correlation term]{Control of $\widetilde Co$.}
We should first study the convergence of the rest term $T_2$ but the mechanisms which intervene in the majoration appear more clearly in $\widetilde Co$. Write $\widetilde Co=\|\widetilde f_{n,\bbeta}-f_0\|_\infty \sum_{j=1}^p \widetilde Co_j(\bbeta_j-\bbeta_{0,j})$ with $\widetilde Co_j=\sum_{i=1}^n \ddot b(\eta_{0,i})X_{i,j}\frac {(\widetilde f_{n,\bbeta}(t_i)-f_0(t_i))}{\|\widetilde f_{n,\bbeta}-f_0\|_\infty}$.

\vspace{10pt}

{\bf Without assumption $(A_{corr})$}

Applying Cauchy-Schwarz leads to $\|\widetilde Co_j\|\leq\|\widetilde f_{n,\bbeta}-f_0\|_n \left(\frac 1 n \sum_i \ddot b(\eta_{0,i}) X_{i,j}^2\right)^{1/2}.$ Under assumptions (A1), (A2) and what precedes, it comes $v_n\|\widetilde Co\|=\bigcirc_\P\left(\sqrt{n}\|\bbeta-\bbeta_0\|\right)$.

\vspace{10pt}

{\bf With $(A_{corr})$}

\vspace{10pt}

We assume the penalty does not deal with scale components of the nonparametric part of the model. As explained in Section A, we consider actually that for all $i=1,\dots,n$, the covariate $\bX_i$ has a null scale representation. Assumption $(A_{corr})$ implies moreover that the polynomial functions $g_j$ have a null wavelet representation and so
$\widetilde Co=\sum_{i=1}^n \ddot b(\eta_{0,i})\left(\bxi_i^{WT}(\bbeta-\bbeta_0)\right)\left(\widetilde f_{n,\bbeta}(t_i)-f_0(t_i)\right)$, where $\bxi_i^W$ represents the wavelets part of $\bxi_i$.

The wavelet transform is orthonormal and thus $\bxi_i^W$ has the same properties as $\xi_i$. When these variables satisfy an sub-gaussian assumption such as (A4.1) or (A4.2), we can apply a control of the term using entropy similar to what has been done before, using \cite{VanderGeer00}.

\begin{itemize}

\item[$\bullet$] When the entropy satisfies $\mathcal H(\mathcal A,\delta,\|.\|_n)\leq A \delta^{-\nu}$, following \cite{MammenVanDerGeer} we get:
$$\sup_{\sqrt{n}\|\bbeta-\bbeta_0\|\leq c} \frac{\widetilde Co_j}{\sqrt{n}\left(\frac{\|\widetilde f_{n,\bbeta}-f_0\|_n}{\|\widetilde f_{n,\bbeta}-f_0\|_\infty}\vee n^{-1/(2+\nu)}\right)^{1-\nu/2}}\,=\,\bigcirc_\P(1).$$
Using $\|\widetilde f_{n,\bbeta}-f_0\|_n=\bigcirc_{\P}(v_n^{-1})$ for $\sqrt{n}\|\bbeta-\bbeta_0\|\leq c$, it comes $\widetilde Co=\bigcirc_\P(1)\left(v_n^{-1}\vee n^{-1/(2+\nu)}\right)^{1-\nu/2}\sqrt{n}\sum_{j=1}^p|\bbeta_j-\bbeta_{0,j}|.$

Using Cauchy-Schwarz inequality,
$$\widetilde Co=\bigcirc_\P(1)\left(v_n^{-1}\vee n^{-1/(2+\nu)}\right)^{1-\nu/2}\sqrt{n}p^{1/2}\|\bbeta-\bbeta_0\|.$$
Provided $\varrho_n=\left(v_n^{-1}\vee n^{-1/(2+\nu)}\right)^{1-\nu/2}p^{1/2}\to 0$, we have $\widetilde Co=o_\P(1)$ for $\bbeta$ such that $\sqrt{n}\|\bbeta-\bbeta_0\|\leq c$.

\item[$\bullet$] When the subgaussian assumption (A4.2) is weakened in the exponential tails assumption (A4.1) we can obtain a similar result using Corallary  8.3 of \cite{VanderGeer00}.

With the $\ell^1$-penalty the majoration for the $\delta$-entropy becomes $\mathcal H(\mathcal A,\delta,\|.\|_n)\leq A \delta^{-2}(\log(n)+\log(1/\delta))$. Note $\tilde\theta_i^W$ the wavelet coefficients of $\tilde f_{n\bbeta}$. Using again \cite{VanderGeer00} we have $$ \P\left(\frac{1}{\sqrt{n\,\log(n)}} \frac{\lambda v_n^2}{\sqrt{n}} \widetilde Co_j\,\geq t \frac{\lambda v_n^2}{n}\sum_i|\tilde\theta_i^W-\theta_{0,i}^W|\right)$$ tends to 0 when $n$ goes to infinity for all $t>C(K_j,\sigma_{0,j}^2)$ with $K_j$ and $\sigma_{0,j}^2$ depending on the distribution of variables $\bxi_{i,j}^{W}$. It was established previously that $\frac{\lambda v_n^2}{n}\sum_i|\tilde\theta_i^W-\theta_{0,i}^W|$ is finite, so we conclude that: $$\widetilde Co=\bigcirc_\P\left(\frac{\sqrt{n}}{v_n^2}p^{1/2}\sqrt{n}\|\bbeta-\bbeta_0\|\right).$$
It is sufficient to suppose $\varrho_n^2=\frac{n}{v_n^4}p\to 0$ to get $\widetilde Co=o_\P\left(\sqrt{n}\|\bbeta-\bbeta_0\|\right)$.

\end{itemize}

\subsection[Control of the rest term]{Control of $T_2$}

The rest term in the Taylor expansion $T_2$ is bounded by: $T_2\leq cst \sum \left|\left(\bX_i^T(\bbeta-\bbeta_0) + \widetilde f_{n,\bbeta}(t_i)-f_0(t_i)\right)^3-\left(\widetilde f_{n,\bbeta}(t_i)-f_0(t_i)\right)^3\right|.$

We decompose in three terms:
\begin{description}

\item[$\bullet$] $T_{21}=\sum \left|\bX_i^T(\bbeta-\bbeta_0)\right|^2|\widetilde f_{n,\bbeta}(t_i)-f_0(t_i)|$ is bounded by
$T_{21}\leq\sup_i\|\bX_i^T\|\|\bbeta-\bbeta_0\|\sum|\bX_i^T(\bbeta-\bbeta_0)||\widetilde f_{n,\bbeta}(t_i)-f_0(t_i)|.$

Proceeding as for the term $\widetilde Co$, when $\sqrt{n}\|\bbeta-\bbeta_0\|\leq c$, we obtain: $T_{21}\leq\sup_i\|\bX_i^T\|\|\bbeta-\bbeta_0\|\bigcirc_\P(\varrho_n).$ 

Assumption (A2) implies $T_{21}=o_\P(\varrho_n).$

\vspace{10pt}

\item[$\bullet$] $T_{22}=\sum \left|\bX_i^T(\bbeta-\bbeta_0)\right|\widetilde f_{n,\bbeta}(t_i)-f_0(t_i)|^2$.
Note that we have $T_{22}\leq\|\widetilde f_{n,\bbeta}-f_0\|_\infty\sum|\bX_i^T(\bbeta-\bbeta_0)||\widetilde f_{n,\bbeta}(t_i)-f_0(t_i)|.$

 And when $\sqrt{n}\|\bbeta-\bbeta_0\|\leq c$, the asymptotic behaviour is $T_{22}=\bigcirc_\P(\|\widetilde f_{n,\bbeta}-f_0\|_\infty\varrho_n).$

\vspace{10pt}

\item[$\bullet$] $T_{23}=\sum \left|\bX_i^T(\bbeta-\bbeta_0)\right|^3.$
Under assumption (A1) and (A2), $T_{23}\leq o(1) h^{1/2}p\left(\sqrt{n} \|\bbeta-\bbeta_0\|_n\right)^3.$

\end{description}

\vspace{10pt}

We have proved that for all $\bbeta$ satisfying $\sqrt{n}\|\bbeta-\bbeta_0\|_n\leq c_2$, the term $T_2$ goes to 0 when $n$ goes to infinity, provided that $\varrho_n\to 0$.

\subsection{Terms with no presence of the functional part}

\begin{itemize}
\item[$\bullet$]{Control of  $B_2=1/2\sum_{i=1}^n \ddot b(\eta_{0,i})\left(\bX_i^T(\bbeta-\bbeta_0)\right)^2.$ }\\
Recall that assumption (A2) stands that $\sum_{i=1}^n \ddot b(\eta_{0,i})\bX_i^T\bX_i$ goes to a strictly positive matrix $\Sigma_0$ when $n$ goes to infinity. Thus, up to a constant $c$, we have $B_2\geq c\gamma(\Sigma_0) n \|\bbeta-\bbeta_0\|^2$, with $\gamma(\Sigma_0)$ smallest eigenvalue of $\Sigma_0$.
Consequently,~~ $\inf_{\|\bbeta-\bbeta_0\|=c'_n }B_2\geq c\gamma(\Sigma_0) n {c'_n}^2$.

\vspace{10pt}

\item[$\bullet$]{Control of $B_1=\sum_{i=1}^n \dot{\mathcal L}(y_i,\eta_{0,i})\left(\bX_i^T(\bbeta-\bbeta_0)\right)$ }

$B_1$ is centered and $\E B_1^2=\sum \ddot b(\eta_{0,i}) \left(X_i^T (\beta-\beta_0)\right)^2.$ Consequently,~~ $\,\sup_{\sqrt{n}\|\bbeta-\bbeta_0\|=c'_n} B_1=\bigcirc_\P(ph^{1/2}c'_n).$

\end{itemize}

\section{Conclusion}

We have proved that $R_n(\bbeta):=\left(K_{n,\lambda}( \widetilde f_{n,\bbeta},\bbeta)-K_{n,\lambda}( \widetilde f_{n,\bbeta},\bbeta_0)\right)-B_1+B_2$ goes to 0 for all $\bbeta$ such that $\sqrt{n}\|\bbeta-\bbeta_0\|\leq c$. Using convexity (see \cite{Rockafellar} as before) the convergence is available for the supremum $\sup_{\sqrt{n}\|\bbeta-\bbeta_0\|\leq c}R_n(\bbeta)\to 0.$

Suppose now that $\sqrt{n}\|\hat \bbeta-\bbeta_0\|\geq c_n$ with $c_n\to \infty$. Then we can build a sequence $c'_n$ satisfying $c'_n\to\infty$, $c'_n\leq c_n$ and $\sup_{\sqrt{n}\|\bbeta -\bbeta_0\|_n= c_n'} \, R_n(\bbeta)\to 0.$ Note that $$R_n(\bbeta)=\left(K_{n,\lambda}( \widetilde f_{n,\bbeta},\bbeta)-K_{n,\lambda}( f_0,\bbeta_0)\right)-\left(K_{n,\lambda}( \widetilde f_{n,\bbeta},\bbeta_0)-K_{n,\lambda}(f_0,\bbeta_0)\right)-B_1+B_2.$$ We have $K_{n,\lambda}( \widetilde f_{n,\bbeta},\bbeta_0)-K_{n,\lambda}(f_0,\bbeta_0)<0$. Using moreover the studies of terms $B_1$ and $B_2$, we obtain that $$\sup_{\sqrt{n}\|\bbeta -\bbeta_0\|_n= c_n'} \, \left(K_{n,\lambda}( \widetilde f_{n,\bbeta},\bbeta)-K_{n,\lambda}(  \widetilde f_{n,\bbeta},\bbeta_0)\right)\,\leq\, \bigcirc_\P(c'_n)-k {c'_n}^2,$$ with $k$ strictly positive constant, and thus~~ $\sup_{\sqrt{n}\|\bbeta -\bbeta_0\|_n= c_n'} \, \left(K_{n,\lambda}(\widetilde f_{n,\bbeta},\bbeta)-K_{n,\lambda}( \widetilde f_{n,\bbeta},\bbeta_0)\right)\,<\, 0$ ~~with a probability going to 1 when $n$ goes to infinity.

Convexity gives: $$\sup_{\sqrt{n}\|\bbeta -\bbeta_0\|_n\geq c_n'} \, K_{n,\lambda}(\widetilde f_{n,\bbeta},\bbeta)-K_{n,\lambda}( \widetilde f_{n,\bbeta},\bbeta_0)<0$$ with a probability going to 1. As the estimator $\hat\bbeta_n$ minimizes $ K_{n,\lambda}(\widetilde f_{n,\bbeta},\bbeta)$, we conclude $$\P(\sqrt{n}\|\hat \bbeta_n -\bbeta_0\|\geq c_n)\to 0.$$


\end{document}